\documentclass[review]{elsarticle}
\usepackage{amsmath, amssymb, amsthm, mathabx, mathrsfs}

\usepackage{lineno,hyperref}

\usepackage{graphicx}

\modulolinenumbers[5]

\journal{Journal of \LaTeX\ Templates}

\newtheorem{theo}{Theorem}[section]

\newtheorem{pro}{Proposition}[section]

\newtheorem{lemma}{Lemma}[section]

\newtheorem{defi}{Definition}[section]

\newtheorem{conj}{Conjecture}[section]

\newtheorem{remark}{Remark}[section]

\numberwithin{equation}{section}

\begin{document}

\begin{frontmatter}

\title{Fredholm Theory for Pseudoholomorphic Curves with Brake Symmetry}

\author[address]{Beijia Zhou}
\ead{beijiachow@gmail.com}

\author[address]{Chaofeng Zhu\corref{mycorrespondingauthor}}
\cortext[mycorrespondingauthor]{Corresponding author}
\ead{zhucf@nankai.edu.cn}

\address[address]{Chern Institute of Mathematics, Nankai University and LPMC, Tianjin 300071, China.}

\begin{abstract}
In this paper we study pseudoholomorphic curves with brake symmetry in symplectization of a closed contact manifold. We introduce the pseudoholomorphic curve with brake symmetry and the corresponding moduli space. Then we get the virtual dimension of the moduli space.
\end{abstract}

\end{frontmatter}

\section{Introduction}

A brake orbit is a periodic orbit of a Hamiltonian system with brake symmetry, which have been studied by mathematicians for decades. In 1948, Seifert has found a brake orbit in a region $G = \{ x | U(x) \leq E \}$, which is homeomorphic to a unit ball, in the condition that the potential energy $U(x_1, \ldots, x_n)$ is an analytic function in a region of ${\bf R}^n$, $\nabla U$ is never vanishing on the bounadry of $G$, and the kinetic energy $T(x, \dot{x}) = \sum_{i,j = 1}^{n}a_{ij}\dot{x}_i\dot{x}_j$ is a positive definite symmetric quadratic form with analytic coefficients \cite{Seifert}. In 2006, Long, Zhang and Zhu\cite{LZZ} used Maslov index to study brake orbit.

Pseudoholomorphic curves in symplectic manifolds have been introduced by Gromov\cite{Gromov} in 1985, and Hofer later in collaboration with Wysocki and Zehnder used pseudoholomorphic curves in symplectization of contact manifolds to study closed Reeb orbits in a series of papers\cite{Hofer},\cite{Hofer1},\cite{Hofer2},\cite{Hofer3}. Similar to  Floer homology in symplectic manifolds, contact homology can be defined in contact manifolds by counting pseudoholomorphic curves. Pseudoholomorphic curves can also be used to study the brake orbit, we refer to the works of Frauenfelder and Kang\cite{Urs} and J.Kim, S.Kim and Kwon\cite{Kim}.

Following these ideas, we hope that contact homology of brake orbits could be defined by counting pseudoholomorphic curves with brake symmetry in the future. In this paper we introduce the pseudoholomorphic curve with brake symmetry and the corresponding moduli space. The main result of this paper is the computation of the virtual dimension of these moduli spaces. It is the first step to define contact homology for brake orbits.

Let $q;q'_1, \ldots, q'_{s^-}$ be brake orbits and  $p'_1, \ldots, p'_{t^-}$ be closed Reeb orbits, see the definition\ref{defi1} and the following paragraph. We abbreviate the muduli space of genus $0$ pseudoholomorphic curves with brake symmetry with one positive puncture $u$, which is asymptotic to a brake orbit $q$, negative punctures $u'_1, \ldots, u'_{s^-}$, which are asymptotic to brake orbits $q'_1, \ldots, q'_{s^-}$, and negative punctures $v'_1, \ldots, v'_{t^-}, w'_1, \ldots, w'_{t^-}$, which are asymptotic to pairs of periodic Reeb orbit $p'_1(t), \ldots, p'_{t^-}(t), Np'_1(-t), \ldots, Np'_{t^-}(-t)$ by $\mathcal{M}_{1+s^-, t^-}(q; q'_1, \ldots, q'_{s^-};\emptyset; p'_1,\\ \ldots, p'_{t^-})$. We define the virtual dimension of $\mathcal{M}=$  the Fredholm index of $D'_F(F \in {\mathcal M})$ $+$ dimension of Teichm\"{u}ller space $\mathcal{T}$ $-$ dimension of automorphism group of Riemann surface with punctures $\Theta$.
\begin{theo}
	The virtual dimension of  moduli space $\mathcal{M}_{1+s^-, t^-}(q; q'_1, \ldots, q'_{s^-};\\ \emptyset; p'_1, \ldots, p'_{t^-})$ is $\frac{n-3}{2} + \mu_1(q) - \sum_{i'=1}^{s^-}\left( \frac{n-3}{2} + \mu_1(q'_{i'}) \right) -\sum_{k'=1}^{t^-}\left( (n-3) + \mu_{CZ}(q'_{k'}) \right).$
\end{theo}
\begin{remark}
The moduli space $\mathcal M$ is not necessarily compact but can be compactified due to SFT compactness. See references \cite{Compact} , \cite{EGH}.
\end{remark}

The precise meaning of the notations in the theorem will be explained in \S \ref{section3}.

The theorem above is a special case of the main result theorem\ref{e} in \S \ref{section3}. There we actually allow more than just one positive puncture and positive genus as well. However, for the applications we have in mind the statement of Theorem 1.1 is sufficient.

\section{Background}

Let
$J_0 = \left(
\begin{array}{cc}
0 & -I_n \\
I_n & 0 \\
\end{array}
\right)$
and
$N_0 = \left(
\begin{array}{cc}
-I_n & 0 \\
0 & I_n \\
\end{array}
\right)$
, where $I_n$ is the identity matrix on $\mathbf{R}^n$. For $H \in C^2(\mathbf{R^{2n}}\backslash \{0\}, \mathbf{R}) \cap C^1(\mathbf{R^{2n}}, \mathbf{R})$ satisfying
\begin{equation}
H(N_0 x) = H(x), \quad \forall x \in \mathbf{R^{2n}},
\end{equation}

we consider the following Hamiltonian system problem for all $t \in \mathbf{R}:$
\begin{align}
\begin{cases}\label{100}
\dot x(t) & = JH'(x(t)), \\
H(x(t)) & = h, \\
x(-t) & = N_0x(t), \\
x(\tau + t) & = x(t).
\end{cases}
\end{align}

A solution $(\tau, x)$ of equations (\ref{100}) is called a \textbf{brake orbit} on the hypersurface $\Sigma := \{ y \in \mathbf{R}^{2n} | H(y) = h\}$. Two brake orbits $(\tau_i, x_i), i=1, 2$, are equivalent if the two brake orbits are geometrically the same, i.e., $x_1(\mathbf{R}) = x_2(\mathbf{R})$. We denote by $[(\tau, x)]$ the equivalence class in this equivalence relation.

Let $(\tau, x)$ be a solution of equations (\ref{100}), We consider the following boundary value problem of the linearized Hamiltonian system at $x$:
\begin{equation}\label{101}
\dot y(t) = JH''(x(t))y(t)
\end{equation}

Denote by $\gamma_{x}$ the fundamental solution of equation (\ref{101}), i.e.
$\gamma_{x}$ is the solution of the problem:
\begin{align}
\dot \gamma_{x}(t) & = JB(t)\gamma_{x}, \\
\gamma(0) & = I_{2n},
\end{align}
where $B(t) = H''(x(t))$ for all $t \in \mathbf{R}$. Then $\gamma_{x}(t) \in \textrm{Sp}(2n).$

\subsection{The Maslov indices}

There are different ways to define the Maslov index, in the literature of contact homology people usually use the way of Robbin, Salamon\cite{Robbin}, and in the book \cite{Long} Professor Long takes almost the way of Cappell, Lee, and Miller\cite{Cappell}. In this paper we will use the first way.

In papers \cite{LZZ}, \cite{LZ} and book \cite{Long}, the authors give the definitions of indices $(i, \nu), (\mu_1, \nu_1), (\mu_2, \nu_2)$ for a symplectic path. Here we briefly review them.

For any $M \in \textrm{Sp}(2n)$, its graph is defined by
$$
\textrm{Gr}(M)  = \{ (x, Mx)| x \in \bf R^{2n}\}
$$

Let

\begin{align*}
L_1 = \{ 0 \} \times {\bf{R}}^n, L_2 = {\bf{R}}^n \times \{ 0 \}\\
W = \{ (x, x) \in {{\bf{R}}^{4n}}| x \in {\bf R}^{2n} \}
\end{align*}

\begin{defi}\label{index}
For any continuous path $\Phi : [a, b] \rightarrow \textrm{Sp}(2n)$ with $\Phi(a) = I_{2n}$, we define the following Maslov-type indices of $\Phi$ for k = 1,2, by

\begin{eqnarray}
  \mu_k^{CLM}(\Phi, [a, b]) &=& \mu^{CLM}(L_k, \Phi L_k), \\
  \nu_k^{CLM}(\Phi, [a, b]) &=& \dim(\Phi(b) L_k \cap L_k),\\
  i(\Phi, [a, b]) &=& \mu^{CLM}(W, Gr(\Phi)) - n,\\
  \nu(\Phi, [a, b]) &=& \dim\ker (\Phi(b) - Id),\\
  \mu_k^{RS}(\Phi, [a, b]) &=& \mu^{RS}(L_k, \Phi L_k), \\
  \mu_{CZ}^{RS}(\Phi, [a, b]) & = & \mu^{RS}(W, Gr(\Phi)) - n
\end{eqnarray}

where $\mu^{RS} , \mu^{CLM}$ are Maslov indices for a pair of Lagrange pathes which are defined in Robbin, Salamon \cite{Robbin} and Cappell, Lee, and Miller \cite{Cappell} respectively.

In this paper we define the Maslov index for $(\tau, x)$ via its associated symplectic path $\gamma_x$ as follows for $k = 1, 2:$
\begin{eqnarray}
\mu_k(x, [0, \tau]) &=& \mu_k^{RS}(\gamma_x, [0, \frac{\tau}{2}]),\\
\nu_k(x, [0, \tau]) &=& \nu_k(\gamma_x, {[0, \frac{\tau}{2}]}),\\
\mu_{CZ}(x, [0, \tau]) &=& \mu_{CZ}(\gamma_x, [0, \tau]),\\
\nu(x, [0, \tau]) &=& \nu(\gamma_x, [0, \tau]).
\end{eqnarray}
\end{defi}
\begin{remark}
In nondegenerate case we have $i (x) = \mu^{RS}_{CZ} (x)$.
\end{remark}

\subsection{Contact Homology}
Contact manifold is the odd dimension analogy of symplectic manifold. It is defined as the following,

\begin{defi}
Let $\Sigma$ be a compact manifold of dimension $2n-1$ with hyperplane distribution $\xi$ in tangent bundle $T\Sigma$. If there exists a 1-form $\alpha$ on $\Sigma$, such that $\alpha \wedge (d\alpha)^{n-1}$ is a volume form, and $\xi = ker \alpha$. We say that $(\Sigma, \xi)$ is a contact manifold and the $\alpha$ is a contact form for $\xi$.
\end{defi}

\begin{remark}
  Given a contact form $\alpha$ for $\xi$, if we multiply $\alpha$ by an any nonvanishing function on $\Sigma$, $f(x) \neq 0, \forall x \in \Sigma$, then $f\alpha$ is another contact form.
\end{remark}

An almost complex structure $J$ on $\xi$, $J : \xi \to \xi, J^2 = -1$ is said to be compatible with $d\alpha$, if it satisfies the conditions $d\alpha(J\cdot, J\cdot) = d\alpha(\cdot, \cdot)$, and $d\alpha(\cdot, J\cdot) > 0$. The set of compatible complex structure is nonempty and contractible, by the same method in \cite[Propsition 2.6.4]{MS1}. Hence the vector bundle $\xi$ over $\Sigma$ is a symplectic vector bundle with symplectic 2-form $d\alpha$, and $d\alpha(\cdot, J\cdot)$ defines an inner product on $\xi$.

 \begin{defi}\label{defi1}
The Reeb vector field $R_{\alpha}$ that is associated to a contact form $\alpha$ is characterized by
\begin{equation*}
\begin{cases}
\begin{array}{ccc}
  i(R_\alpha)d\alpha & = & 0, \\
   \alpha(R_\alpha) & = & 1.
\end{array}
\end{cases}
\end{equation*}
\end{defi}

An orbit $x(t): {\bf R} \to \Sigma$ which satisfies the condition $\frac{d}{dt}(x(t)) = R_\alpha(x(t))$ is called a Reeb orbit. A closed Reeb orbit is a Reeb orbit which satisfies the condition that there exists $\tau > 0$, such that $x(t+\tau) = x(t)$. A closed Reeb orbit is called a periodic Reeb orbit as well. Let $\phi(t)$ be the flow generated by $R_\alpha$, i.e. $\frac{d\phi(t)}{dt} = R_\alpha \circ \phi$. $\phi(t)$ preserves $\alpha$ and $d\alpha$, $\phi_t^*\alpha = \alpha, \phi_t^*d\alpha = d\alpha$. So $\phi(t)$ preserve $\xi$ and the linearized map $d\phi(t): \xi_{x(0)} \to \xi_{x(t)}$ is a symplectic map. See \cite[page 1634]{Siefring}. We can take a trivialization along the orbit i.e. it is a symplectic isomorphism from $\xi_{x(t)}$ to standard linear symplectic space ${\bf R}^{2n-2}$, $\varphi(t): \xi_{x(t)} \to {\bf R}^{2n-2}$. Hence the map $\gamma_x(t) = \varphi(t) \circ d\phi_t \circ \varphi(0)^{-1} : [0, \tau] \to \textrm{Sp}(2n)$ is a symplectic path with $\gamma_x(0) = I_{2n}$, and we define indices $\mu_{CZ}(x), \mu_1(x), \mu_2(x)$ for $x$ by $\mu_{CZ}^{RS}(\gamma_x), \mu_1^{RS}(\gamma_x), \mu_2^{RS}(\gamma_x)$ in Definition \ref{index}.

If we take a different trivialisation $\varphi'(t)$ along the orbit $x(t)$, then we have a loop in symplectic group $\Upsilon(t) : \varphi' \circ \varphi^{-1}(t) : [0, \tau] \to \textrm{Sp}(2n)$, and we have the symplectic path $\gamma'_x(t) = \varphi'(t) \circ d\phi_t \circ \varphi'(0)^{-1} : [0, \tau] \to \textrm{Sp}(2n)$ associated to the trivialisation $\varphi'(t)$. It is a well know fact that a symplectic loop $\Upsilon(t)$ is homotopic to a standard loop $diag(R(\frac{2\pi kt}{\tau}), 1, \cdots, 1)$, where $R(t) = \left(
                                                                                                                   \begin{array}{cc}
                                                                                                                     \cos t & -\sin t \\
                                                                                                                     \sin t & \cos t \\
                                                                                                                   \end{array}
                                                                                                                 \right)
$,  $k$ is called the degree of $\Upsilon(t)$. It is easy to see that $\mu_{CZ}(\gamma'_x) = \mu_{CZ}(\gamma_x) + 2\textrm{degree}\Upsilon(t)$. See\cite[Remark 5.3]{Robbin}.

In the following we assume each closed Reeb orbit $x(t)$ of $\Sigma$ is nondegenerate, i.e. for each period $\tau$ of $x(t)$($\tau$ does not need to be the minimal period), the symplectic matrix $d\phi_\tau: \xi_{x(0)} \to \xi_{x(\tau)}$ has no eigenvalue $1$.

For every closed Reeb orbit $x(t)$ in $\Sigma^{2n-1}$, we define $|x| = \mu_{CZ}(\gamma_x) + n-3$.
Let $x(t)$ be be a closed Reeb orbit with minimal period $\tau$. We consider the multiple covers $x^m$ of $x$ with period $m\tau, m \geq 2$. There are two ways that the grading of   $x^m$ can behave :
\begin{enumerate}[(i)]
\item the parity of $|x^m|$ is the same for all $m \geq 1.$
\item the parity for the even multiples $|x^{2k}|, k\geq 1$, disagrees with the parity for the odd multiples $|x^{2k-1}|, k\geq 1$.
\end{enumerate}

 \begin{defi}
In the second case, the even multiples $x^{2k}, k \in N^+$, are called bad orbits. An orbit that is not bad is called good.
\end{defi}

\begin{defi}
Let $C_*$ be the graded $\mathbf Z$ module, which is freely generated by all good closed Reeb orbits $x$ and $|x| = *$.
\end{defi}

 \begin{defi}
 The symplectization of the contact manifold $(\Sigma, \xi)$ with contact form $\alpha$ is the symplectic manifold $({\bf R} \times \Sigma, d(e^t\alpha))$, where $t$ is the coordinate of $\bf R$.
 \end{defi}

We can extend the almost complex structure $J$ on $\xi$ to an almost structure on ${\bf R} \times \Sigma$ (we still denote by $J$) compatible with $\omega = d(e^t\alpha))$ by defining $J\frac{\partial}{\partial t} = R_\alpha$.

The contact homology counts the number of pseudoholomorphic curves in the symplectization.
Let $\Theta = S^2\backslash \{u, v_1, \ldots, v_s\},$ where $u, v_1, \ldots, v_s$ are distinct points on $S^2$. Let $j$ be a complex structure on $\Theta$. Then $(\Theta, j)$ is a Riemann surface with punctures. Let map $F : (\Theta, j) \rightarrow ({\bf R}\times \Sigma, J)$ satisfies $dF \circ j = J \circ dF$. We write $F = (a, f)$, where $a$ is the component in ${\bf R}$, $f$ is the component in $\Sigma$. Let $(\rho, \theta)$ be polar coordinates centered on a puncture. On a small neighbourhood of the puncture, we require
\begin{align*}
\lim_{\rho\to 0}a(\rho, \theta) & = +\infty \quad \text{for puncture}\ u,\\
& = -\infty \quad \text{for punctures}\  v_1, \ldots, v_s, \\
\lim_{\rho\to 0}f(\rho, \theta) & = x(-\frac{T}{2\pi}\theta) \quad \text{for puncture}\ u, \\
& = x_i(\frac{T}{2\pi}\theta) \quad \text{for punctures}\ v_1, \ldots, v_s,
\end{align*}
for some closed orbits $x$ of period $T$ and $x_i$ of period $T_i, i = 1, \ldots, s.$ Here $u$ is called positive singularity and $u_1, \ldots, u_s$ are called negative singularity. We define the energy $E(F) := \textrm{sup}\{\int_\Theta F^*d(\phi\alpha)| \phi : {\bf R} \to [0, 1], \phi' > 0 \}$, the finite energy condition $0 < E(F) < \infty$ for a pseudoholomorphic curve guarantees the above requirements, see \cite[Theorem 1.2]{Hofer1}.

Two holomorphic curves $F : S^2\backslash \{u, v_1, \ldots, v_s\} \to (\mathbf{R} \times M, J)$ and $F' : S^2\backslash \{u', v'_1, \ldots, v'_s\} \to (\mathbf{R} \times M, J)$ are equivalent if and only if there exists a biholomorphism $h : (S^2, j) \to (S^2, j')$ so that $h(u) = u', h(v_i) = v'_i$ for $i = 1, \ldots, s$, and $F = F' \circ h$.

\begin{defi}
The moduli space of pseudoholomorphic curves $\mathcal{M}(x; x_1, \ldots, x_s)$ is the set of equivalence classes of pseudoholomorphic curves with above asymptotic conditions. It has an $\mathbf{R}$-action induced by the translation $t \to t + \Delta t$ in $\mathbf{R} \times M$.
\end{defi}

According to results of thesis of Bourgeois \cite[Corollary 5.4]{Bourgeois} and result of Dragnev \cite[Theorem 1.8]{Dragnev}, the dimension of moduli space is shown by the following theorem
\begin{pro}
The dimension of moduli space $\mathcal{M}(x; x_1, \ldots, x_s)/{\bf R}$ is
\begin{equation*}
(n - 3)(1 - s) + \mu_{CZ}(x) - \sum_{i=1}^s \mu_{CZ}(x_i) - 1 + 2 c_1(\Sigma, F^*\xi)
\end{equation*}
Where $c_1(\Sigma, F^*\xi)$ is the first Chern number of $F^*\xi$ on $\Sigma$, relative to given trivialization of $F^*\xi$ along the closed Reeb orbits at the punctures.
\end{pro}

\begin{remark}
Taking a different trivialization, the Maslov index $\mu_{CZ}$ and $c_1$ can change but the dimension of Moduli space shown by the formula will not change.
\end{remark}

The rational weights take into account the automorphisms of pseudoholomorphic curves : if $F$ is an element in $\mathcal{M}/\mathbf{R}$ of dimension $0$, then the weight of $F$, $weight(F)$ is $\frac{1}{k}$, where $k$ is the order of the automorphism group of $F$.

The number
\begin{equation*}
  n_{x_a, x_b} = \left\{
     \begin{array}{ll}
       0 & \hbox{if dim}\mathcal{M}(x_a; x_b) \neq 1 \\
       \sum_{F\in\mathcal{M}(x_a; x_b)/\bf{R}}weight(F) & \hbox{if dim} \mathcal{M}(x_a; x_b) = 1
     \end{array}
   \right.
\end{equation*}
counts pseudoholomorphic cylinders joining closed orbit $x_a$ and $x_b$.

We difine
\begin{gather*}
 d: C_k \to C_{k-1} \quad \text{for all}\ k \in \mathbf{Z} \\
  dx_a = m_a \sum_{b, |b| = |a| - 1} n_{x_a, x_b}x_b
\end{gather*}
where $m_a$ is the multiplicity of $\gamma_a$.

\begin{conj}[Eliashberg, Hofer\cite{EGH} ]

If $C_k = 0$ for $k = -1, 0, 1$, then
\begin{enumerate}[(i)]
  \item $d^2 = 0$,
  \item $H_*(C_*, d)$ is independent of the contact form $\alpha$ for $\xi$, and the complex structure $J$.
\end{enumerate}
\end{conj}

\begin{defi}
The cylindrical contact homology $HC_*(\Sigma, \xi)$ is the homology $H_*(C_*, d)$ of the chain complex $(C_*, d)$.
\end{defi}

The following facts connect the periodic orbit of Hamiltionian system and the periodic Reeb orbit in contact manifold.

 \begin{lemma}
The solutions of Hamiltonian system $\dot z(t) = JH'(z(t))$ in ${\bf R}^{2n}$ are equivalent to the Reeb orbits of contact form $\alpha = \sum_{i = 1}^{n}\frac{1}{2}(-y_i dx_i + x_i dy_i)$ on the energy surface $\Sigma = \{H^{-1}(1) \}$, provided that $x_iH_{x_i} + y_iH_{y_i} \neq 0.$
\end{lemma}
\begin{proof}
We only need to prove the Hamiltonian vector field is proportional to the Reeb vector field.
On one hand $X_H = (-H_y, H_x)$. Here $H_x = (H_{x_1}, \ldots, H_{x_n})$, $H_y = (H_{y_1}, \ldots, H_{y_n})$. On the other hand $d\alpha = \sum^n_{i=1}dx_i \wedge dy_i$, $i(-H_y, H_x)d\alpha = -H_{y_i}dy_i -H_{x_i}dx_i $, which equals $0$ on $T\Sigma$. $\alpha(-H_y, H_x) = \frac{1}{2}(y_iH_{y_i} + x_iH_{x_i}) \neq 0$. Hence $R_\alpha = fX_H, f \in C(M, {\bf R})$.
\end{proof}

The brake orbits are closed characteristics with brake symmetry $x(-t) = Nx(t)$. Since the Hamiltonian vector field is proportional to the Reeb vector field, the brake orbits are equivalent to the closed Reeb orbits with brake symmetry.

We will show that the hamiltonian vector field and the Reeb vector field of $\alpha = \sum_{i = 1}^{n}\frac{1}{2}(-y_i dx_i + x_i dy_i)$ on the surface $\Sigma, N\Sigma = \Sigma$ have a symmetry:
\begin{lemma}\label{a1}
If $H(N_0x) = H(x)$, then $X_H(N_0x) = -N_0 X_H(x)$.
If $N^*\alpha = -\alpha$, then $R_\alpha |_{Nx} = -NR_\alpha |_x$.
\end{lemma}
\begin{proof}
If the hamiltonian function satisfies $H(N_0 x) = H(x)$, then $N_0 H'(N_0 x) = H'(x)$,  $-N_0 JH'(N_0 x) = JH'(x)$, $JH'(N_0 x) = -N_0 JH'(x)$. So we get $X_H(N_0 x) = -N_0 X_H(x)$.

If we assume that $N^*\alpha = -\alpha$, we have $i(-NR_\alpha)_{Nx}d\alpha_{Nx} = i(N_* (-R_\alpha)_{x})d\alpha_{Nx}\\ = i((-R_\alpha)_{x})(N^*d\alpha)_x = i((-R_\alpha)_{x})(-d\alpha)_x = 0$, and $i(-NR_\alpha)_{Nx}\alpha_{Nx}\\ = i(N_* (-R_\alpha)_{x})\alpha_{Nx} = i((-R_\alpha)_{x})(N^* \alpha)_x = i((-R_\alpha)_{x})(-\alpha)_x = 1$. So we get $R_\alpha |_{Nx}  = -NR_\alpha |_x $.
\end{proof}

If we deform a energy surface $\Sigma \subset {\bf R}^{2n}, N_0 \Sigma = \Sigma$, in the meanwhile we keep the brake symmetry of it $N_0 \Sigma_t = \Sigma_t$, the Hamiltonian vector field and the Reeb vector field of $\alpha = \sum_{i = 1}^{n}\frac{1}{2}(-y_i dx_i + x_i dy_i)$ on $\Sigma_t$ will keep the relation $X_{N_0 x} = -N_0 X_x$. So a brake orbit on an energy surface $\Sigma \subset {\bf R}^{2n}, N_0 \Sigma = \Sigma$ is equivalent to a Reeb orbit with brake symmetry for some contact form with involution $(\alpha, N)$.

Let $\alpha$ be a contact form on $\Sigma$, $N$ is an involution on $\Sigma$, they satisfy $N^*\alpha = -\alpha$.
\begin{align}
  \dot x(t) & = R_\alpha(x), \label{102}\\
  x(-t) & = Nx(t), \\
  x(\tau + t) & = x(t). \label{103}
\end{align}

We also call a solution $(\tau, x)$ of equations (\ref{102}) - (\ref{103}) a brake orbit.

\begin{remark}
Define $L = \{z \in \Sigma|Nz = z\}$. By the uniqueness of solutions to ODEs, the requirement of equations (\ref{102}) - (\ref{103}) are also equivalent to
\end{remark}

\begin{align}
	\dot x(t)  = R_\alpha(x), \label{104} \\
	x(0), x(\frac{\tau}{2}) \subset L. \label{105}
\end{align}

A fact is we cannot choose the start point $x(0)$ arbitrarily along a brake orbit, this is different to a periodic orbit. Because $x(0)$ is a fix point of involution $N$. By lemma \ref{a1} we know $R_\alpha |_{Nx} = -NR_\alpha |_x$, we get $R_\alpha = -NR_\alpha$ along $L$. Hence $R_\alpha$ transverse to the tangent space of $L$, $TL = \{X \in T\Sigma | X = NX\}$. Therefore we get a lemma, which will be used later

\begin{lemma}\label{a2}
Because $R_\alpha$ is transverse to $TL$ along $L$, a brake orbit is not invariant under $S^1$ action $x(t) \to x(t+\Delta t)$.
\end{lemma}

\section{Main Result}\label{section3}
Following the construction in contact homology, we will give the construction of muduli space of pseudoholomorphic curves with brake symmetry.

Let $(\Sigma, \alpha)$ be a contact manifold with the contact structure $\xi = ker \alpha$. Let $N$ be an involution of $\Sigma$, $N^2 = \textrm{Id}$, such that $N^*\alpha = -\alpha$. We consider brake orbits on $(\Sigma, \alpha, N)$. We assume all periodic Reeb orbits are nondegenerate.

\begin{remark}
For a brake orbit $q$, it is a brake orbit and a periodic orbit, there are two kind of nondegeneracy $\nu(q) = 0$ or $\nu_1(q)= 0$ for $q$. We have the relation $\nu_1(q) + \nu_2(q) = \nu(q),$ which can be proved as the proof of Proposition C in \cite{LZZ}. In this paper, they are all equal to $0$, according to our nondegeneracy assumption.
\end{remark}

We take the symplectization of $\Sigma$, $(M = {\bf R} \times \Sigma, d(e^t\alpha))$, where $t$ is the coordinate of $\bf R$..We take an almost complex structure $J$ on $\xi$ such that $NJ = -JN$ and compatible with $\alpha$, then we extend $J$ to $M$ by letting $J\frac{\partial}{\partial t} = R_\alpha$ which is compatible with $d(e^t\alpha)$. We also extend the involution $N$ in $\Sigma$ to ${\bf R} \times \Sigma$ by $\textrm{Id} \times N$. For convenience, we still use $N$ for the extension $\textrm{Id} \times N$. Note that $L = \{z \in \Sigma|Nz = z\}$ is a Legendrian submanifold of $\Sigma$.

Using the method of contact homology, we consider the pseudoholomorphic curves as follows.

Let $\Xi_g$ be a Riemann surface of genus $g \geq 0$ with an involution $\mathcal{N}$. Where $\mathcal{N} : \Xi_g \to \Xi_g$, $\mathcal{N}^2 = \mathrm{Id}$. We abbreviate $\Xi_g$ by $\Xi$ if we do not need to mention the genus. We denote $\mathcal{L} := \{p \in \Xi_g | \mathcal{N}p = p \}$ to be the fixed point set of $\mathcal{N}$. Let $m$ be the number of connected component of $\mathcal{L}$. Let $\langle\Xi\rangle = 1$, if the quotient $\Xi/\mathcal{N}$ is orientable, i.e. $\Xi\backslash \mathcal{L}$ has two connected components, and $\langle\Xi\rangle = 0$ otherwise. For example, let $\Xi_0$ be a sphere in ${\bf R}^3$, and define $\mathcal{N}: S^2 \to S^2$ by $\mathcal{N}(x, y,z) = (-x, y, z)$, where $(x, y, z)$ is the coordinate of $p$, $\forall p \in S^2$. Two Riemann surfaces with involution are equivalent if there is a biholmorphism $\tau$ between them such that $ \tau\mathcal{N}_1 = \mathcal{N}_2\tau$. The equivalent class of Riemann surface with involution is characterised by $(g, m, \langle\Xi\rangle)$. For fix genus $g$, there are $[\frac{g+1}{2}]$ equivalent class with $\langle\Xi\rangle = 1$ and $g+1$ equivalent class with $\langle\Xi\rangle = 0$. See \cite[Corollary 1.1]{Natanzon}. We consider only the Riemann surface that has $m \geq 1$ and $\langle\Xi\rangle = 1$ in the following.

Let $\{u_1, \ldots, u_{s^+}; u'_1, \ldots, u'_{s^-}; v_1, \ldots, v_{t^+}, w_1, \ldots, w_{t^+}; v'_1, \ldots, v'_{t^-}, w'_1, \ldots, w'_{t^-}\}$ be punctures in $\Xi_g$ such that $\mathcal{N}(u_i) = u_i, \mathcal{N}(u'_{i'}) = u'_{i'}$, for $i = 1, \ldots, s^+, i' = 1, \ldots, s^-$, $\mathcal{N}(v_k) = w_k, \mathcal{N}(v'_{k'}) = w'_{k'}$, for $k = 1, \ldots, t^+, k' = 1, \ldots, t^-$, namely $u_i, u'_i$ are fixed point of $\mathcal{N}$, $\mathcal{N}$ permutes $v_k$ and $w_k$, and $\mathcal{N}$ permutes $v'_{k'}$ and $w'_{k'}$. We take $j$ a complex structure of $\Xi_g$, such that $\mathcal{N}j = -j\mathcal{N}$. We denote $\Theta_{g, s^+ +s^-, t^+ +t^-} := \Xi_g\backslash\{u_1, \ldots, u_{s^+}; u'_1, \ldots, u'_{s^-}; v_1, \ldots, v_{t^+}, w_1, \ldots, w_{t^+}; v'_1, \ldots,\\ v'_{t^-}, w'_1, \ldots, w'_{t^-}\}$ be the Riemann surface with $s^+ +s^-$ symmetric punctures and $t^+ +t^-$ symmetric pairs of punctures. Without loss of generality, we can assume the punctures $v_1, \ldots, v_{t^+}, v'_1, \ldots, v'_{t^-},$ are in one component of $\Theta_g \backslash \mathcal{L}$ and $w_1, \ldots, w_{t^+}, w'_1, \ldots, w'_{t^-}$ are in the other component. We abbreviate $\Theta_{g, s^+ +s^-, t^+ +t^-}$ by $\Theta$ when we need not to mention genus and number of punctures.

We consider a map $F = (a, f) : (\Theta, j, \mathcal{N}) \to ({\bf R} \times \Sigma, J, N)$, such that $dF \circ j = J \circ dF$, $a(\mathcal{N}\cdot) = a(\cdot)$ and $f(\mathcal{N}\cdot) = Nf(\cdot)$, where $a$ is the component in ${\bf R}$, $f$ is the component in $\Sigma$. Let $(\rho, \varphi)$ be the polar coordinates centered on a puncture, we want
\begin{equation}\label{a}
\begin{aligned}
  \lim_{\rho\to 0}a(\rho, \varphi) & = +\infty \quad \text{for punctures}\ u_1, \ldots, u_{s^+}, v_1, \ldots, v_{t^+}, w_1, \ldots, w_{t^+}\\
   & = -\infty \quad \text{for punctures}\  u'_1, \ldots, u'_{s^-}, v'_1, \ldots, v'_{t^-}, w'_1, \ldots, w'_{t^-} \\
  \lim_{\rho\to 0}f(\rho, \varphi) & = q_i(-\frac{T_i}{2\pi}\varphi) \quad \text{for punctures}\ u_1, \ldots, u_{s^+} \\
   & = q'_{i'}(\frac{T'_{i'}}{2\pi}\varphi) \quad \text{for punctures}\ u'_1, \ldots, u'_{s^-} \\
   & = p_k(-\frac{T_k}{2\pi}\varphi) \quad \text{for punctures}\ v_1, \ldots, v_{t^+} \\
   & = Np_k(\frac{T_k}{2\pi}\varphi) \quad \text{for punctures}\ w_1, \ldots, w_{t^+} \\
   & = p_k'(\frac{T'_{k'}}{2\pi}\varphi) \quad \text{for punctures}\ v'_1, \ldots, v'_{t^-} \\
   & = Np'_{k'}(-\frac{T'_{k'}}{2\pi}\varphi) \quad \text{for punctures}\ w'_1, \ldots, w'_{t^-}
\end{aligned}
\end{equation}
for some brake orbits $q_i(t)$ of period $T_i, i = 1, \ldots, s^+$,  $q'_{i'}(t)$ of period $T'_{i'}, i' = 1, \ldots, s^-$ and closed Reeb orbits $p_j(t), Np_j(-t)$ of period $T_{s^+ + k},  k = 1, \ldots, t^+$, $p'_{k'}(t), Np'_{j'}(-t)$ of period $T'_{s^- + k'},  k' = 1, \ldots, t^-$. We call such $F$ pseudoholomorphic curve with brake symmetry. $u_1, \ldots, u_{s^+}, v_1, \ldots, v_{t^+}, w_1, \ldots, w_{t^+}$ are called positive punctures, and $u'_1, \ldots, u'_{s^-}, v'_1, \ldots, v'_{t^-}, w'_1, \ldots, w'_{t^-}$ are called negative punctures.

We also introduce the equivalent relations.
Two pseudoholomorphic curves with brake symmetry $F : (\Theta, j, \mathcal{N}) \to (M, J, N)$, $\Theta = \Xi_g \backslash \{u_1, \ldots, u_{s^+}; u'_1, \ldots, u'_{s^-};\\ v_1, \ldots, v_{t^+}, w_1, \ldots, w_t; v'_1, \ldots, v'_{t^-}, w'_1, \ldots, w'_{t^-}\}$ and $F' : (\widetilde\Theta, \tilde j,\widetilde{ \mathcal{N}}) \to (M, J, N)$, $\widetilde\Theta = \widetilde\Xi_g \backslash \{\tilde{u}_1, \ldots, \tilde{u}_{s^+}; \tilde{u}'_1, \ldots, \tilde{u}'_{s^-}; \tilde{v}_1, \ldots, \tilde{v}_{t^+}, \tilde{w}_1, \ldots, \tilde{w}_{t^+}; \tilde{v}'_1, \ldots, \tilde{v}'_{t^-}, \tilde{w}'_1, \ldots, \tilde{w}'_{t^-}\}$\\ are equivalent if and only if there exists a biholomorphism with brake symmetry $h : (\Theta, j, \mathcal{N}) \to (\widetilde\Theta, \tilde j,\widetilde{ \mathcal{N}}), h(\mathcal{N}\cdot) = \mathcal{N}h(\cdot)$ so that $h(u_i) = \tilde{u}_i$, for $i = 1, \ldots, s^+$, $h(u'_{i'}) = \tilde{u}'_{i'}$, for $i' = 1, \ldots, s^-$, $h(v_{s + k}) = \tilde{v}_{s + k}$, for $k = 1, \ldots, t^+$, $h(v'_{s' + k'}) = \tilde{v}'_{s' + k'}$, for $k' = 1, \ldots, t^-$ and $F = F' \circ h$.

\begin{defi}\label{b}
The moduli space of pseudoholomorphic curves of genus $g$ with brake symmetry $\mathcal{M}_{g, s^+ +s^-, t^+ +t^-}(q_1, \ldots, q_{s^+}; q'_1, \ldots, q'_{s^-}; p_1, \ldots, p_{t^+}; p'_1, \ldots, p'_{t^-})$ is the set of equivalent classes of pseudoholomorphic curves of genus $g$ with brake symmetry which converges to brake orbits and closed Reeb orbits as (\ref{a}). It has an $\mathbf{R}$-action induced by the translation $t \to t + \Delta t$ in $\mathbf{R} \times \Sigma$. When $g = 0$, we abbreviate it by $\mathcal{M}_{s^+ +s^-, t^+ +t^-}$.
\end{defi}

\subsection{Asymptotic Behavior}
In this subsection, $F = (a, f): (\Theta, j, \mathcal{N}) \to (M, J, N)$ is a pseudoholomorphic curve with brake symmetry such that $f$ is asymptotic to a brake orbit near a neighborhood of a positive puncture.

We first give the asymptotic behavior of pseudoholomorphic half cylinder with brake symmetry such that $f$ is asymptotic to a brake orbit. Our method is similar to \cite{Hofer1}, \cite{Mora}, \cite{Siefring}.

For a closed Reeb orbit $p(t)$ of period $\tau$, from the standard method in contact manifold, see \cite[Lemma 2.3]{Hofer1} and \cite[Theorem 2.5.15]{Geiges}, there exists a neighborhood $U$ of $S^1 \times {\bf R}^{2n-2}$ and a neighborhood $V$ of the periodic orbit $p$ in $\Sigma$, such that there exists a diffeomorphism $\psi: U \to V, g: U \to \bf R$,  with $\psi^*\alpha = g\alpha_0$, where $\alpha_0 = d\theta + \sum_{i=1}^{n} x_idy_i$, $\theta$ is the coordinate of $S^1$, $(x_1, \cdots, x_{n-1}, y_1, \cdots, y_{n-1})$ is the coordinate of ${\bf R}^{2n-2}$, and $g$ satisfies $g|(\theta, 0, 0) = \tau, dg|(\theta, 0, 0) = 0$. Without loss of generality, we can assume $\tau = 1$.

We take the polar coordinate near a neighborhood of a positive puncture $u$,
\begin{align*}
\varphi: [R, +\infty) \times S^1 \to {\bf C}\backslash(0, 0) \\ (s, t) \mapsto e^{-2\pi(s+it)}
\end{align*}
where we set $u$ as $(0, 0)$ in the complex plane $\bf C$.

Taking $R$ sufficiently large, $F$ can be represented by these coordinates
\begin{gather*}
F_1 = (a \circ \varphi, \psi^{-1} \circ F \circ \varphi) : [R, +\infty) \times S^1 \to {\bf R} \times S^1 \times {\bf R}^{2n-2},\\
(s, t) \mapsto \\ (a(s,t), \theta(s,t), x(s, t)_1, \cdots, x(s, t)_{n-1}, y(s, t)_1, \cdots, y(s ,t)_{n-1})
\end{gather*}

The Cauchy-Riemann condition $dF \circ j = J \circ dF$ can be written as
\begin{align}
\begin{cases}
a_s - \alpha(f_t) = 0,\\
a_t + \alpha(f_s) = 0,\\
\pi(f_s) + J(f)\pi(f_t) = 0. \label{1}
\end{cases}
\end{align}
where $\pi$ is the projection to $\xi$.

We write $x = (x(s, t)_1, \cdots, x(s, t)_{n-1})$, $y = (y(s, t)_1, \cdots, y(s ,t)_{n-1}))$. We define $R_\alpha := (X_1, X, Y)$ in this coordinate, where $X$ is the $x$ component, $Y$ is the $y$ component of $R_\alpha$ respectively. The third equation of  (\ref{1}) is equivalent to
\begin{equation}
  \left(
     \begin{array}{c}
       x_{s} \\
       y_{s} \\
     \end{array}
   \right) + J(s, t)\left(
     \begin{array}{c}
       x_{t} \\
       y_{t} \\
     \end{array}
   \right) + a_t\left(
        \begin{array}{c}
          X \\
          Y \\
        \end{array}
      \right) - a_s J(s, t) \left(
           \begin{array}{c}
             X \\
             Y \\
           \end{array}
         \right) = 0.
         \label{2}
\end{equation}
Where $J(s,t)$ is the matrix which represent the almost complex structure $J$ in this coordinate.

We define
\begin{equation*}
  Z(t, x, y) = \left(
     \begin{array}{c}
       X \\
       Y \\
     \end{array}
   \right) = D(t, x, y) \left(
                          \begin{array}{c}
                            x \\
                            y \\
                          \end{array}
                        \right)
\end{equation*}
with
\begin{equation*}
  D(t, x, y) = \int_0^1 dR_\alpha(t, \eta x,\eta y)d\eta
\end{equation*}

We write $z = \left(
\begin{array}{c}
x \\
y
\end{array}
\right)$, then the equation (\ref{2}) can be written as
\begin{equation}
z_{s} + J(s, t) z_{t} + (a_t - a_s J(s, t))D(t,x, y) z = 0\label{d}.
\end{equation}

By \cite[Equation (20)]{Hofer1} and \cite[Equation (136) (137)]{Mora}, we can find the $(a, \theta)$ component satisfies the equation
\begin{align}\label{6}
\begin{cases}
a_s = (\theta_t + x \cdot y_t)g \\
a_t = -(\theta_s + x \cdot y_s)g
\end{cases}
\end{align}
where $x \cdot y$ means the standard Euclidean product of two vectors.

We can always reduce the almost complex structure $J(s, t)$ into the standard matrix $J_0$. See \cite[Remark 2.9]{Hofer1}. We define $S(s, t) := (a_t - a_s J_0)D(t, x, y)$, and the operator $A(s) := -J_0\frac{d}{dt} - S(s, t) : W^{1, 2}(S^1, \xi) \to L^2(S^1, \xi)$. $A(s)$ has the limit $A_\infty = -J_0 \frac{d}{dt} - J_0 dR_\alpha$, where $S_\infty = - J_0 dR_\alpha$ is an $S^1$ orbit of symmetric matrices, using the same discussion in \cite[page354-355]{Hofer1}. By a similar proof as \cite[Lemma 2.5]{Hofer1}, the dimension of $\textrm{ker}A_\infty$ is equal to the multiplicity of eigenvalue $1$ for symplectic matrix $\Phi(1)$, where $\Phi(t)$ is the fundamental solution of equation $\dot{y}(t) = J_0 S_\infty(t)y$.

If $p$ is a brake orbit, $R_\alpha$ has brake symmetry $R_\alpha(\mathcal{N}\cdot) = -NR_\alpha(\cdot)$. We get $NdR_\alpha(\mathcal{N}\cdot)N = -dR_\alpha(\cdot)$, and $NJdR_\alpha(\mathcal{N}\cdot)N = JdR_\alpha(\cdot)$, i.e. $N_0 S_\infty(-t) N_0 = S_\infty(t)$, where by taking coordinate appropriately the involution $N$ on $\xi$ along $p$ can be represented by the standard matrix
$$N_0 = \left(
\begin{array}{cc}
-I_n & 0 \\
0 & I_n \\
\end{array}
\right)$$

We define a similar operator for a brake orbit $\widetilde A(s) := -J_0 \frac{d}{dt} - S(s, t) : {\widetilde W}^{1, 2}(S^1, \xi) \to {\widetilde L}^2(S^1, \xi),$ where ${\widetilde W}^{1, 2}(S^1, \xi) := \{w\in W^{1, 2}(S^1, \xi)|  w(-t) = Nw(t)\}$, ${\widetilde L}^{2}(S^1, \xi) := \{\ell \in L^{2}(S^1, \xi)|  \ell(-t) = N\ell(t)\}$. We have a similar result, that  the dimension of ker ${\widetilde A}_\infty$ is equal to the intersection of two space $\Phi(b) L_1$ and $L_1$, namely dim$(\Phi(b) L_1 \cap L_1)$, where $\Phi(t)$ is the fundamental solution of equation $\dot{y}(t) = JS_\infty(t)y$, $L_1 = 0 \times {\bf R}^n$.

By equation (\ref{6}) and $f$ is asymptotic to $p$, the $(a, \theta)$ component satisfies the Cauchy-Riemann equation $(\frac{\partial}{\partial s} + J \frac{\partial}{\partial t})(a, \theta) = 0$ in the limit.

If $F = (a, f)$ satisfies the finite energy condition $E(F) < \infty$, $f$ is asymptotic to some periodic orbit. It is a well known fact that $A$ is a selfadjoint operator. Under the nondegenerate condition $\nu(p) = 0$, $A$ does not have eigenvalue $0$. By result of Hofer\cite{Hofer1}, Mora \cite{Mora}, and Siefring \cite{Siefring}, we have:
\begin{pro}\label{3}
	\begin{equation*}
	z(s,t) = e^{\lambda s}(e(t) + r(s,t) )
	\end{equation*}
	Where $\lambda$ is a negative eigenvalue of $A_\infty$, $e(t) \neq 0$ is an eigenvector with eigenvalue $\lambda$, and $\lim_{s \to 0} r(s, t) = 0$.

There exists $d > 0, s_0, \theta_0, a_0$ such that, for every multi-index $I$ there is a constant $c_I$ so that
	
	\begin{eqnarray*}
		|\partial^I r(s, t)| & \leq & c_I e^{-ds},  \\
		|\partial^I (\theta(s, t) - t - \theta_0)| & \leq & c_I e^{-ds}, \\
		|\partial^I (a(s, t) -\tau s - a_0)| & \leq & c_I e^{-ds}
	\end{eqnarray*}
	for all $s \geq s_0$.
\end{pro}

Therefore, we say
\begin{defi}\label{convergent}
	The function $F_1   : [R, +\infty) \times S^1 \to {\bf R} \times S^1 \times {\bf R}^{2n-2}$, $F_1 = (a(s,t), \theta(s,t), x(s, t)_1, \cdots, x(s, t)_{n-1}, y(s, t)_1, \cdots, y(s ,t)_{n-1})$ is $(d, k, p)$-convergent to the periodic orbit $p$, if there exist $a_0, \theta_0$, such that $e^{ds}|a(s, t) -\tau s - a_0|$, $e^{ds}|\theta(s, t) - t - \theta_0|$, $e^{ds}z \in W^{k,p}$.

A function $F: \Theta \to M$ is called $(d, k, p)$-convergent to brake orbits $q_1, \ldots, q_s;\\ q'_1, \ldots, q'_{s'}$ and closed Reeb orbits $p_1, \ldots, p_t; p'_1, \ldots, p'_{t'}$, if $F \in W^{k, p}(\Theta, M)$, $F$ is convergent to the periodic orbits as in equation (\ref{a}), and at the punctures $F$ is $(k, p, d)$-convergent to the corresponding orbit. We denote the space of $(k, p, d)$-convergent functions by $W^{k, p, d}(\Theta, M).$
\end{defi}

These are the results for pseudoholomorphic curves $F$ which are asymptotic to given closed Reeb orbits. If we consider pseudoholomorphic curves with brake symmetry $F(\mathcal{N}\cdot) = NF(\cdot)$ which  are asymptotic to given brake orbits. We still have the result of Proposition \ref{3} for $F$.

We consider the map $d\mathcal{N} : T\Theta \to T\Theta$. By $\mathcal{N}^2 = Id$ we have splitting $T\Theta = T\Theta^+ \oplus T  \Theta^-$, where $T\Theta^{\pm}$ is the eigenspace of $\pm 1$ for $d\mathcal{N}$. Similarly by $N^2 = 1$ we have $TM = TM^+ \oplus TM^-$, where $TM^{\pm}$ is the eigenspace of $\pm 1$ for $dN$. The relation $dF \circ d\mathcal{N} = dN \circ dF$ implies $dF: T\Theta^{\pm} \to TM^{\pm}$.

\subsection{Fredholm Theory}

In this subsection we calculate the virtual dimension of moduli space of pseudoholomorphic curves with brake symmetry.

We follow the method of Bourgeois\cite{Bourgeois} and Dragnev\cite{Dragnev}.

Let $q_1, \ldots, q_{s^+}; q'_1, \ldots, q'_{s^-}$ be brake orbits and  $p_1, \ldots, p_{t^+}; p'_1, \ldots, p'_{t^-}$ be closed Reeb orbits. We have introduced the moduli space $\mathcal{M}_{g, s^+ +s^-, t^+ +t^-}(q_1, \ldots,\\ q_{s^+}; q'_1, \ldots, q'_{s^-}; p_1, \ldots, p_{t^+}; p'_1, \ldots, p'_{t^-})$ in Definition \ref{b} . Then each $F \in \mathcal{M}$ satisfies the description in Proposition \ref{3}.

Let $\bar{\mathcal{B}} = \bar{\mathcal{B}}^{k, p, d}(q_1, \ldots, q_{s^+}; q'_1, \ldots, q'_{s^-}; p_1, \ldots, p_{t^+}; p'_1, \ldots, p'_{t^-}),d > 0, p > 2$, $k$ is positive integer, be the Banach manifold of maps $F = (a, f): \Theta \to \mathbf{R} \times \Sigma$ with brake symmetry $F(\mathcal{N}\cdot) = NF(\cdot)$, and $F$ is $(k, p, d)$-convergent to the correspond closed orbit near each puncture.

We follow the standard pseudoholomorphic curves theory explained by McDuff and Salamon\cite{MS}.
For each $F \in \mathcal{M}$, we denote the Banach space $\{G \in W^{k-1, p, d}(\Theta, \Lambda^{0, 1} \otimes F^*TM)| G(\mathcal{N}\cdot) = NG(\cdot) \}$, the $(k-1, p, d)$-convergent $(0, 1)$-form on $\Theta$ with value in $F^*TM$ and brake symmetry, by $\bar{W}^{k-1, p, d}(\Theta, \Lambda^{0, 1} \otimes F^*TM)$.
We define $\bar{\mathcal{E}}$ the Banach bundle over $\bar{\mathcal{B}}$, the fiber $\bar{\mathcal{E}}_{F}$ over $F$ will be the Banach space $\bar{\mathcal{E}}_F = \bar{W}^{k-1, p, d}(\Theta, \Lambda^{0, 1} \otimes F^*TM)$.

The nonlinear operator $\bar{\partial} : \bar{\mathcal{B}} \to \bar{\mathcal{E}}, \bar{\partial}F = dF + J \circ dF \circ j$ is a section of this bundle. Its zero set is the space of pseudoholomorphic curves with brake symmetry.

We allow the conformal  structure of Riemann surface to vary, so we consider the enlarged Banach manifold $\widetilde{\mathcal{B}} = \mathcal{T}_{g, s^+ +s^-,t^+ +t^-} \times \bar{\mathcal{B}}$, where $\mathcal{T}_{g, s^+ +s^-,t^+ +t^-} $ is the Teichm\"{u}ller space for Riemann sphere with brake symmetry $\Theta_{g, s^+ +s^-,t^+ +t^-}$ as we defined above. Similarly we have a bundle $\widetilde{\mathcal{E}}$ over $\widetilde{\mathcal{B}}$ and the Cauchy-Riemann operator induces a section $ \widetilde{\mathcal{B}} \to \widetilde{\mathcal{E}} : (u, j) \to \bar\partial u, \bar{\partial} F = dF + J \circ dF \circ j$.

The zero set $\bar{\partial}^{-1}(0)$ of the Cauchy-Riemann section $\bar{\partial} : \widetilde{\mathcal{B}} \to \widetilde{\mathcal{E}}$ is the set of pseudoholomorphic curves with brake symmetry. Two pseudoholomorphic curves with brake symmetry $F : (\Theta, j) \to (M, J)$ and $F' : (\tilde\Theta, j') \to (M, J)$ are equivalent if and only if there exists a biholomorphism with brake symmetry $h : (\Theta, j) \to (\tilde\Theta, j'), h(N\cdot) = h(\cdot)$ so that $F = F' \circ h$.

It follows by elliptic regularity that holomorphic curves are smooth maps if we take the almost complex structure $J$ smooth, therefore the definition of moduli space is independent of the values of $k, p, d$ as long as they satisfy the requirements that $d > 0, k \geq 1, p > 2$. Therefore from Definition\ref{b} the moduli space $\mathcal{M}_{g, s^+ +s^-,t^+ +t^-}(q_1, \ldots, q_{s^+}; q'_1, \ldots, q'_{s^-}; p_1, \ldots, p_{t^+}; p'_1, \ldots, p'_{t^-})$ consists of the equivalent classes in $\bar{\partial}^{-1}(0) \subset \mathcal{B}.$

By linearizing the section map $\bar{\partial} : \bar{\mathcal{B}} \rightarrow \bar{\mathcal{E}}$ at the point of pseudoholomorphic curve $F \in \mathcal{B}$, we get operator $$d\bar{\partial}_F : T_F \bar{\mathcal{B}} \to T_{(F, 0)}\bar{\mathcal{E}}$$

The space $T_{(F, 0)}\bar{\mathcal{E}}$ has decomposition $T_{(F, 0)}\bar{\mathcal{E}} = T_F \bar{\mathcal{B}} \oplus \bar{\mathcal{E}}_F$. We define $\pi : T_{(F, 0)}\bar{\mathcal{E}} \to \bar{\mathcal{E}}_F$ to be the projection to second component. We denote  $D'_F := \pi \circ d\bar{\partial}_F$. We have $T_F \bar{\mathcal{B}} \cong {\bf R}^{s^+ + s^- +2t^+ +2t^-} \oplus \bar{W}^{d,k,p}(\Theta, F^*TM)$, where $\bar{W}^{d,k,p}(\Theta, F^*TM)$ is the space
 $\{H \in W^{k, p, d}(\Theta, F^*TM)| H(\mathcal{N}\cdot) = NH(\cdot) \},$ the $(k, p, d)$-vector field along $\Theta$ with value in $F^*TM$ and brake symmetry.

 The component ${\bf R}^{s^+ + s^- +2t^+ +2t^-}$ comes from the freedom of image in the direction $\frac{\partial}{\partial s}$ for brake orbit $q_1, \ldots, q_{s^+}, q'_1, \ldots, q'_{s^-}$, and $(\frac{\partial}{\partial s}, R_\alpha)$ for closed Reeb orbit $p_1, \ldots, p_{t^+}, p'_1, \ldots, p'_{t^-}$. A brake orbit has only $\frac{\partial}{\partial s}$ direction, because a brake orbit is not $S^1$ action invariant by lemma \ref{a2}, a brake orbit has only $\frac{\partial}{\partial s}$ direction. Let $\rho(s)$ be a function with support in a small neighborhood of symmetric puncture, depending only on s and equal to $1$ for $s$ sufficient large. Let $\sigma(s)$ be a function with support in a small neighborhood of a symmetric pair of  punctures, depending only on s and equal to $1$ for $s$ sufficient large at one puncture and satisfies $\sigma(\mathcal{N}\cdot) = \sigma(\cdot)$. The component ${\bf R}^{s^+ + s^- +2t^+ +2t^-}$ is spanned by vectors $\rho\frac{\partial}{\partial s}$ at each symmetric puncture and $(\sigma\frac{\partial}{\partial s}, \sigma R_\alpha)$ at each symmetric pair of punctures.

 Note that, because of the exponential behavior of $F$, the linear operator is exponentially converging to its asymptotic value at each puncture. Hence, the image of the summand ${\bf R}^{s^+ + s^- +2t^+ +2t^-}$ in the domain is contained in $\bar{W}^{k-1, p, d}(\Theta, \Lambda^{0, 1} \otimes F^*TM)$.

In the following we first consider the component $D'_F|\bar{W}^{d,k,p}$, we denote $D_F = D'_F|\bar{W}^{d,k,p}$. Without loss of generality, we consider $D_F$ near a positive puncture in the special coordinate $(a(s,t), \theta(s,t), x(s, t)_1, \cdots, x(s, t)_{n-1}, y(s, t)_1, \cdots,\\ y(s ,t)_{n-1})$. From discussion after the equation (\ref{d}) in subsection 3.1, we know that the linearized operator $D_F$ on the component $\xi$ has the form $\frac{\partial}{\partial s} + J_0 \frac{\partial}{\partial t} + S(s, t)$ near a puncture. The operator $A(s)=  -J_0 \frac{\partial}{\partial t} - S(s, t)$ has the limit $A_\infty$. Ker$A_\infty$ is $0$, and $N_0S_\infty(-t)N_0 = S_\infty(t)$. Moreover the $(a, \theta)$ component of $D_F$ near a puncture has limit $\bar\partial = \frac{\partial}{\partial s} + J \frac{\partial}{\partial t}$. So near a puncture, $D_F$ has the limit $D_F^\infty = \bar\partial \oplus \bar\partial + S_\infty(t)$.

Next we do the Morse-Bott construction as in \cite{Dragnev},\cite{Bourgeois}. Let
$\sigma :\\ \bar{W}^{k, p, d}(\Theta, F^*TM) \rightarrow \bar{W}^{k, p}(\Theta, F^*TM)$ and $\sigma' : \bar{W}^{k, p, d}(\Theta, \Lambda^{0, 1}\otimes F^*TM) \rightarrow \bar{W}^{k,p}(\Theta, \Lambda^{0, 1}\otimes F^*TM)$ be the multiplication by $e^{ds}$, which are isomorphism for all $k$. Consider the operator
\begin{align*}
\widetilde{D}_F = \sigma' \circ D_F \circ \sigma^{-1} :\bar{W}^{k,p}(\Theta, F^*TM) \rightarrow \bar{W}^{k-1,p}(\Theta, \Lambda^{0, 1} \otimes F^*TM)
\end{align*}

Near the puncture, we know $D^\infty_F$ has the limit form $\bar\partial \oplus \bar\partial + S_\infty(t)$, so $\widetilde{D}^\infty_F = \sigma' \circ D_F \circ \sigma^{-1}$ has the form $(\bar\partial - d \mathrm{Id}) \oplus (\bar\partial + S_\infty(t) -d \mathrm{Id})$. If we take $d > 0$ sufficient small, $(\bar\partial - d \mathrm{Id})$ and $(\bar\partial + S_\infty(t) -d \mathrm{Id})$ are nondegenerate operators. So $\widetilde{D}^\infty_F$ is a nondegenerate operator.

For the negative puncture, there is a similar result. $\widetilde{D}^{- \infty}_F = \sigma' \circ D_F \circ \sigma^{-1}$ has the form $(\bar\partial + d \mathrm{Id}) \oplus (\bar\partial + S_\infty(t) + d \mathrm{Id})$, where $\sigma, \sigma'$ will be multiplication by $e^{-ds}$.

\begin{lemma}
The operator
\begin{align*}
\widetilde{D}_F = \sigma' \circ D_F \circ \sigma^{-1} :\bar{W}^{k,p}(\Theta, F^*TM) \rightarrow \bar{W}^{k-1,p}(\Theta, \Lambda^{0, 1} \otimes F^*TM)
\end{align*}
is a Fredholm operator.
\end{lemma}

\begin{proof}
We consider the operator
\begin{align*}
\widetilde{D}^{asym}_F = \sigma' \circ D_F \circ \sigma^{-1} :{W}^{k,p}(\Theta, F^*TM) \rightarrow {W}^{k-1,p}(\Theta, \Lambda^{0, 1} \otimes F^*TM)
\end{align*}
In \cite[3.1.11 Theorem]{Schwarz}, M.Schwarz proved the operator $\widetilde{D}^{asym}_F$ of the same form without brake symmetry satisfies the inequality
\begin{equation*}
  ||u||_{W^{1,p}} \leq C(||\widetilde{D}^{asym}_F u||_{L^P} + ||K u||_{Z})
\end{equation*}
for all $u \in W^{1, p}$, some constant $C > 0$ and a compact operator $K : W^{1, p} \to Z$ with some Banach space $Z$.

Since $\bar{W}^{1, p} \subset W^{1, p}$ and $\widetilde{D}_F$ has the same form as $\widetilde{D}^{asym}_F$, the same inequality holds for $\widetilde{D}_F$.

Then by \cite[3.1.10 lemma]{Schwarz}, we know $\widetilde{D}_F$ is semi-Fredholm, i.e. its kernel is finite-dimensional and its range is closed. The adjoint of $\widetilde{D}_F$ has the same form as $\widetilde{D}_F$, the adjoint of $\widetilde{D}_F$ is a semi-Fredholm operator as well. Therefore we get $\widetilde{D}_F$ is a Fredholm operator.
\end{proof}

Next we calculate the index of $\widetilde{D}^\infty_F$, by using the pair of pants induction in Schwarz's thesis\cite{Schwarz} and the Appendix C in McDuff and Salamon\cite{MS}.

We first give a definition,
\begin{defi}
A positive cap $D^+$ is defined by gluing the positive half cylinder to the unit disk,
\begin{align*}
   D^+ = D^2 & \cup_f Z^+, D^2 = \{|z| \leq 1\}, Z^+ = [0, \infty) \times S^1 \subset {\bf C}/i{\bf Z},\\
  & f : {\bf C}/i \to {\bf C}\backslash{0}, (s, t) \mapsto e^{2\pi(s+it)} = re^{i\varphi}.
\end{align*}

A negative cap $D^-$ is defined by gluing the negative half cylinder to the unit disk,
\begin{align*}
   D^- = D^2 & \cup_f Z^+, D^2 = \{|z| \leq 1\}, Z^- = (-\infty, 0] \times S^1 \subset {\bf C}/i{\bf Z},\\
  & \bar f : {\bf C}/i{\bf Z} \to {\bf C}\backslash{0}, (s, t) \mapsto e^{2\pi(-s+it)} = re^{i\varphi}.
\end{align*}
\end{defi}

We briefly review Schwarz's procedures for calculating the index. For details, see \cite[Chapter 3]{Schwarz}. Let $E$ be a complex bundle of rank $n$ over a Riemann surface with punctures $\Sigma\backslash\{ v_1, \cdots, v_s \}$, $D$ is a Cauchy-Riemann operator over $\Sigma\backslash\{ v_1, \cdots, v_s \}$. Let $S(t)$ be a loop of $2n \times 2n$ symmetric matrices, and $\Phi(t) : [0, 1] \to \rm{Sp}(2n)$ be the fundamental solution of the equation $\dot y = S(t)y$. We define $\mu_{CZ}(S(t)) := \mu_{CZ}(\Phi(t))$.

Step 1, we calculate the Cauchy-Riemann operator $D = \bar{\partial} + \beta(s)ds \otimes \omega \mathrm{Id}$ over a positive cap, where $\beta$ is define in the below \ref{bb}. And we get the index $D = n + \mu_{CZ}(\omega \mathrm{Id})$ for positive cap and index $D = n - \mu_{CZ}(\omega \mathrm{Id})$ for a negative cap.

Step 2, if the Riemann surface $\Sigma$ is glued by two Riemann surface $\Sigma_{*}, * =1, 2,$ along a common closed boundary $L$, i.e. $\Sigma = \Sigma_1 \bigcup_L \Sigma_2$, the Cauchy-Riemann operator $D$ on $\Sigma$ has the same form along $L$ as $D$ on $D|\Sigma_{*}, * = 1,2$, which is represented by $A= -J_0\frac{\partial}{\partial t} - S(t)$. Then we have index ${D}|\Sigma =$ index ${D}|\Sigma_1 + $ index ${D}|\Sigma_2$.

 Step 3, using this gluing formula, we can glue a positive cap in first step by a cylinder with positive boundary $-J_0\frac{\partial}{\partial t} - S(t)$ and negative boundary $-J_0\frac{\partial}{\partial t} - \omega \textrm{Id}$. We know the index of $D$ on cylinder is the difference $\mu_{CZ}(S(t)) - \mu_{CZ}(\omega Id)$ by spectrum flow. So we get the index of $D$ with general boundary condition $-J_0\frac{\partial}{\partial t} - S(t)$ on a positive cap is $n + \mu_{CZ}(S(t))$. Similarly the index on a negative cap is $n-\mu_{CZ}(S(t))$.

 Step 4, from Riemann-Roch Theorem, we know index of $D$ on a compact Riemann surface $\Sigma$ is $n\chi(\Sigma) + 2c_1(E)$.

 Step 5, we can glue the Riemann surface $\Sigma\backslash\{ v_1, \cdots, v_s \}$ with punctures by opposite caps of the same boundary condition. From the gluing formula, we get index $D = n(\chi(\Sigma) - s) + \sum_i \mu_{CZ}(v^+_i) - \sum_j \mu_{CZ}(v^-_j)$.

Our first step is the index theorem for a positive cap. Let $(E, N)$ be a complex vector bundle with brake symmetry. We can assume a trivialisation $T: E \to D^+ \times {\bf C}^n$ of $E$ satisfies $T(\bar z) = \bar T(z)$ by our assumption and taking appropriate coordinate.

Let $S(t)$ be a loop of $2n \times 2n$ symmetric matrices with brake symmetry $N_0S(-t)N_0 = S(t)$, and $\Phi(t) : [0, 1] \to \rm{Sp}(2n)$ be the fundamental solution of the equation $\dot y = S(t)y$. We define $\mu_1(S(t)) := \mu_1(\Phi(t))$, $\nu_1(S(t)) := \nu_1(\Phi(t))$. We give the first Lemma.

\begin{lemma}
	Let $E$ be a complex bundle of rank $n$ over a positive cap $D^+$ with brake symmetry, and we use our trivialisation of $E$ by $T$. $\bar{\partial} = \frac{\partial}{\partial s} + J \frac{\partial}{\partial t}$ be the standard Cauchy-Riemann operator. We define
\begin{align}\label{bb}
\beta(s) = \begin{cases}
0, & s \leq 2, \\
1 & s \geq 3,
\end{cases} \text{and strictly monotone on}\ (2, 3)
\end{align}

We define the operator $\bar{D} = \bar{\partial} + \beta(s)ds \otimes\omega \mathrm{Id}: \bar {W}^{k, p}(D^+, E) \to \bar{W}^{k-1, p}(D^+,\\ \Lambda^{0, 1} \otimes E)$, where  $\bar{W}$ is the symmetirc sections of $E$ with brake symmetry $u(\bar z) = \bar u(z)$, for $u \in \bar {W}$ and we require $\frac{\omega}{2\pi} \notin {\bf Z}$. Then the Fredholm index of $\bar{D}$ is $\frac{n}{2} + \mu_1(\omega \mathrm{Id})$.
\end{lemma}

\begin{proof}
	Because $E$ can split into $n$ rank $1$ complex vector bundles, we can consider only the rank $1$ complex bundle. And from the bootstrap procedures of elliptic operator, the index does not depend on $k, p$ if $kp > 2.$ We can reduce the theorem to the case $k = 1, p > 2$. So in the following we consider $\bar{D} : \bar {W}^{1, p}(D^+, E) \to \bar{W}^{0, p}(D^+, \Lambda^{0, 1} \otimes E)$, where $p > 2$.
	
In the discussion in \cite[section 3.3]{Schwarz}, $\bar{D}$ is represented by $\frac{\partial}{\partial s} + i\frac{\partial}{\partial t} + \beta\cdot \omega \mathrm{Id}$. We can represent any section of $E$ by its Fourier series $u = A_k(r)e^{ik\phi}$, where $(r, \theta)$ is the polar coordinate. From the proof of \cite{Schwarz} proposition of section 3.3.5, we have a bijection between the kernel and cokernel of the operator $D = \bar{\partial} + \omega \mathrm{Id}$ without our brake symmetry condition and a complex linear space ${\bf C}^h$, where $h$ is some positive integer.

If $u \in \text{ker}D$, then $u$ satisfies the conditions $A_k = c_k r^k, \text{when}\ 0 < r < \frac{1}{2},$ and $A_k = c_k r^{(k - \frac{\omega}{2\pi})}, \text{when}\ r >3$. A necessary and sufficient condition for $u \in \text{ker}D$ is $k \geq 0$ and $2k\pi < \omega$.
\begin{equation*}
\text{dim}_{\bf C} \text{ker}D = \text{cardinality of} \{k \in {\bf Z}| 0 \leq k <  \frac{\omega}{2\pi}\}
\end{equation*}
And if $\omega < 0$, ker$D$ is $\{0\}$.

If $u \in \text{coker}D$, then $u$ satisfies the conditions $A_k = c_k r^{-k-1}, \text{when}\ 0 < r < \frac{1}{2},$ and $A_k = c_k r^{ \frac{\omega}{2\pi} - k}, \text{when}\ r >3$. A necessary and sufficient condition for $u \in \text{coker}D$ is $k \leq -1$ and $\omega < 2k\pi$.
\begin{equation*}
\text{dim}_{\bf C} \text{coker}D = \text{cardinality of} \{k \in {\bf Z}| \frac{\omega}{2\pi} < k \leq  -1\}
\end{equation*}
And if $\omega  > -1$, coker$D$ is $\{0\}$.

After the same procedure and adding the brake symmetry condition, $u(\bar z) = \bar u(z)$, we get in our case $c_k \in {\bf R}$.
	
Set $l-1 = [\frac{\omega}{2\pi}]$, then index $\bar{D} = l $.
	
	Since the boundary condition of the operator $\bar{D}$ is $A_\infty = -J\frac{d}{dt} - \omega \mathrm{Id}$, the correspondent symplectic path for $\bar{D}$ is $\Phi(t) = R(\omega t) = \left(
                                                                                           \begin{array}{cc}
                                                                                             \cos(\omega t) & -\sin(\omega t) \\
                                                                                             \sin(\omega t) & \cos(\omega t) \\
                                                                                           \end{array}
                                                                                         \right)
, t \in[0, 1]$, then the Maslov index $\mu_1(\omega \mathrm{Id})$ is Maslov index $\mu_1$ for symplectic path $\Phi(t)= R(\omega t), t \in[0, \frac{1}{2}]$. Using the result of \cite[Example 3.2]{LZZ}, we know $\mu_1(\omega \mathrm{Id}) = \frac{1}{2} + [\frac{\omega}{2}\cdot \frac{1}{\pi}] = \frac{1}{2} + l -1$.\footnote{Here we use the Robbin, Salamon index, it is a little different to definition in \cite{LZZ}, we take just the half the crossing dimension at the start point, the other steps is the same as \cite{LZZ}}
	
Eventually, we get the equality $l = \text{index} \bar{D} = \frac{1}{2} + \mu_1(\omega Id)$
	
\end{proof}

Let $(\Theta_{*}, \mathcal{N}_{*}), * = 1,2$ be two Riemann surface with brake symmetry and punctures, where $\Theta_1 \in \mathcal{M}_{g, s^+_1 +s^-_1,t^+_1 +t^-_1}$, $\Theta_2 \in \mathcal{M}_{g, s^+_2 +s^-_2,t^+_2 +t^-_2}$, and $s^-_1 = s^+_2$, $t^-_1 = t^+_2$. Let $(\Theta, \mathcal{N})$ be the Riemann surface which is glued by $\Theta_1$ and $\Theta_2$ along the negative punctures of $\Theta_1$ and positive punctures of $\Theta_2$,  i.e. $\Theta = \Theta_1 \bigcup_\mathcal{L} \Theta_2$, where $\mathcal{L}$ is the common boundary of $\Theta_1$ and $\Theta_2$, $\mathcal{L}$ has also brake symmetry induced from Riemann surface. Let $\bar{D}$ be a Cauchy-Riemann operator $\bar{D} : \bar{W}^{k, p}(\Theta, E) \to \bar{W}^{k-1, p}(\Theta, \Lambda^{0, 1} \otimes E)$  on $\Theta$, and $\bar{D}|\Theta_{*}, * = 1,2$ on $\Theta_{*}, * = 1,2$. $\bar{D}$ and $\bar{D}_{*}, * = 1,2$ has the same form in a neighborhood of $\mathcal{L}$, which is represented by $A = -J_0\frac{d}{dt} - S(t)$. By almost the same procedure as in \cite[\S 3.2]{Schwarz}, we can get our second lemma.

\begin{lemma}[\textbf{Gluing Formula}]
Given Riemann surfaces and Cauchy-Riemann operators with brake symmetry described in the above, we have
$$\text{index} \bar{D}|\Theta = \text{index} \bar{D}|\Theta_1 +  \text{index} \bar{D}|\Theta_2$$
\end{lemma}

We consider a Cauchy-Riemann operator on a cylinder with brake symmetry
\begin{align*}
\bar{D} & : \bar{W}^{k, p}(S^1 \times \mathbf{R}, E) \to \bar{W}^{k-1, p}(S^1 \times \mathbf{R}, \Lambda^{0, 1} \otimes E)\\
\bar{D} & = \frac{\partial}{\partial s} - {\widetilde A}(s)
\end{align*}
where ${\widetilde A}(s) = -J\frac{\partial}{\partial t} - S(s, t):  \bar{W}^{k, p}(S^1, E) \to \bar{W}^{k-1, p}(S^1, \Lambda^{0, 1} \otimes E)$. At the same time ${\widetilde A}(s)$ has positive and negative infinite limits ${\widetilde A}(\infty) = -J\frac{\partial}{\partial t} - S^\infty(t)$, $A(-\infty) = -J\frac{\partial}{\partial t} - S^{-\infty}(t)$, $S^\infty(t), S^{-\infty}(t)$ are symmetric matrces with brake symmetry $NS^*(-t)N = S^*(t), *= \infty, -\infty$. We require the kernel of the operators ${\widetilde A}(\infty)$ and ${\widetilde A}(-\infty)$ are $\{0\}$. The index of $\bar{D}$ is the difference $\mu_1(S^\infty(t)) -\mu_1(S^{-\infty}(t))$ by spectral flow. Because the sections have brake symmetry, the dimension of ker ${\widetilde A}(s) = -J\frac{d}{dt} - S^\infty(s, t)$ is equal to $\nu_1(S(t))$, which is explained in \cite[Proposition 4.1]{LZZ}. Taking the explanation of Maslov index by spectral flow in \cite[\S 7]{Cappell}, we can also get index of a Cauchy-Riemann operator on the cylinder by spectral flow. That is our third lemma.

\begin{lemma}
Given a Cauchy-Riemann operator $\bar{D}$ on a cylinder with brake symmetry, which is represented by $\bar{D}  = \frac{\partial}{\partial s} - {\widetilde A}(s)$. ${\widetilde A}(s) = -J\frac{d}{dt} - S^\infty(s, t)$ has positive and negative infinite limit $A(-\infty) = -J\frac{\partial}{\partial t} - S^{-\infty}(t)$, and the kernel of the operators ${\widetilde A}(\infty)$ and ${\widetilde A}(-\infty)$ are $\{0\}$. Then we have
$$\text{index}\ \bar D = \mu_1(S^\infty(t)) -\mu_1(S^{-\infty}(t))$$
\end{lemma}

For a general positive cap boundary condition, which has boundary condition  ${\widetilde A}^\infty = -J\frac{\partial}{\partial t} - S^\infty(t)$, $S^\infty(t)$ is a symmetric matrix with brake symmetry $N_0S^\infty(-t)N_0 = S^\infty(t)$. We require the kernel of the operator ${\widetilde A}^\infty$ is $\{0\}$. We can glue a positive cap with boundary condition ${\widetilde A} = -J\frac{\partial}{\partial t} - \omega \mathrm{Id}$ by a cylinder with positive boundary condition ${\widetilde A} = -J\frac{\partial}{\partial t} - S^\infty(t)$ and negative boundary condition ${\widetilde A} = -J\frac{\partial}{\partial t} - \omega \mathrm{Id}$. By Gluing Formula we get the index of $\bar{D}$ on the positive cap with positive boundary condition${\widetilde A} = -J\frac{\partial}{\partial t} - S^\infty(t)$ is $\frac{n}{2} + \mu_1(S^\infty(t))$.

For a negative cap, we can also get index $\bar{D} = \frac{n}{2} - \mu_1(S^{-\infty(t)})$ by similar procedures.

Here is our next lemma.
\begin{lemma}
Let $E$ be a complex bundle of rank $n$ over a positive cap $D^+$ with brake symmetry, and we use our trivialisation of $E$ by $T$. Let $\bar D$ be a Cauchy-Riemann operator with brake symmetry on $E$ over $D^+$, which has boundary condition that ${\widetilde A}^\infty = -J\frac{\partial}{\partial t} - S^\infty(t)$, and the kernel of the operator ${\widetilde A}^\infty$ is $\{0\}$. Then the Fredholm index of $\bar D$ is $\frac{n}{2} + \mu_1(S^{\infty(t)})$.

Similarly the Fredholm index of $\bar D$ on a negative cap $D^-$ is $\frac{n}{2} - \mu_1(S^{-\infty(t)})$.
\end{lemma}

We can glue a positive cap $D^+$ with brake symmetry by a negative cap $D^-$ with brake symmetry, both caps have the same boundary condition $A(t) = -J\frac{d}{dt} - \omega \mathrm{Id}$. Then we get a sphere with brake symmetry by gluing the above two caps. From the Gluing Formula, index $\bar{D}|{S^2} = $ index $\bar{D}|(D^+)\ +$ index $\bar{D}|(D^-) = n$.

The above formula can be achieved by Riemann-Roch formula as well.\\ Riemann-Roch formula shows real dimension of  index $D$ on a closed Riemann surface $\Theta$ is $n\chi(\Theta) + 2c_1(E)$. We abbreviate the space $W^{k, p}(\Theta, E)$ by $W$ and  $\bar{W}^{k, p}(\Theta, E)$ by $\bar{W}$. We define the space $\bar{W}' = \{u \in W | u({\mathcal N}z) = -Nu(z)\}$. We have the following result.

\begin{lemma}
$W$ has decomposition $W = \bar{W} \oplus \bar{W}'$, and the space $\bar{W}$ is isomorphic to the space $\bar{W}'$ by the multiplying of $J$.
\end{lemma}

\begin{proof}
For any $u \in W$, we have the decomposition $u(z) = \frac{u(z) + Nu(\mathcal{N}z)}{2} + \frac{u(z) - Nu(\mathcal{N}z)}{2}$. We denote $u_1 := \frac{u(z) + Nu(\mathcal{N}z)}{2}, u_2 := \frac{u(z) - Nu(\mathcal{N}z)}{2}$. It is easy to check $u_1 \in \bar{W}$ and $u_2 \in \bar{W}'$.

Let $u_1$ be an element in $W_1$. We have $Ju_1(Nz) = JNu_1(z) = -NJu_1(z) = -N(Ju(z))$, so $Ju_1 \in \bar{W}'$ and vice versa. The space $\bar{W}$ is isomorphic to the space $\bar{W}'$ by the multiplying of $J$.
\end{proof}

A Cauchy-Riemann operator $\bar D$ on a closed Riemann surface is represented by $\bar{D} = \bar{\partial} + S(z),$ where $z$ is the coordinate of the surface. Because index of a Cauchy-Riemann operator is invariant under deformation of $S(z)$, we can choose $S(z)$ commute with almost complex structure $J$, $JS = SJ$. Then $\bar D$ is commute with $J$, $J\bar{D} = \bar{D}J$. By the decomposition $W = \bar{W} \oplus \bar{W}'$, the dimensions of range, ker and coker of operator $\bar{D}$ is half of operator $D$, we get the \textbf{Riemann-Roch formula with brake symmetry}: index $\bar{D} = \frac{1}{2}$ ind$D = \frac{n}{2}\chi(\Theta) + c_1(E)$.

\begin{lemma}[\textbf{Riemann-Roch formula with brake symmetry}]
Given a Cauchy-Riemann operator $\bar D$ with brake symmetry on a Riemann surface with brake symmetry $\Theta$,
$$\text{index}\ \bar{D} = \frac{n}{2}\chi(\Theta) + c_1(E)$$
\end{lemma}

\begin{remark}
Warning: For a Riemann surface with punctures, we do not always have ind$\bar{D} = \frac{1}{2}\text{ind}D$, because the boundary condition $A(t) = -J_0\frac{d}{dt} - S(t)$, we cannot always choose $S$ commute with $J$. For example, a positive cap with boundary condition $-J\frac{d}{dt} - S(t)$, $\mu_1(S(t)) \neq \frac{1}{2}\mu_{CZ}(S(t)).$
\end{remark}

On the other hand, the Riemann-Roch formula with brake symmetry can be recovered by pair of pants reduction. We can easily get the index on a pair of pants. Then we cut $\Theta$ into a union of several pairs of pants with brake symmetry. Using the Gluing Formula, the Riemann-Roch formula with brake symmetry can be recovered by adding the index on each pair of pants together.

For example, we can cut a torus into four parts, which consist of a negative cap, a pair of pant with a positive puncture and two negative puncture,  a pair of pant with two positive puncture and a symmetric puncture and a positive cap.

\begin{figure}
  \centering
  \includegraphics[scale = 0.4]{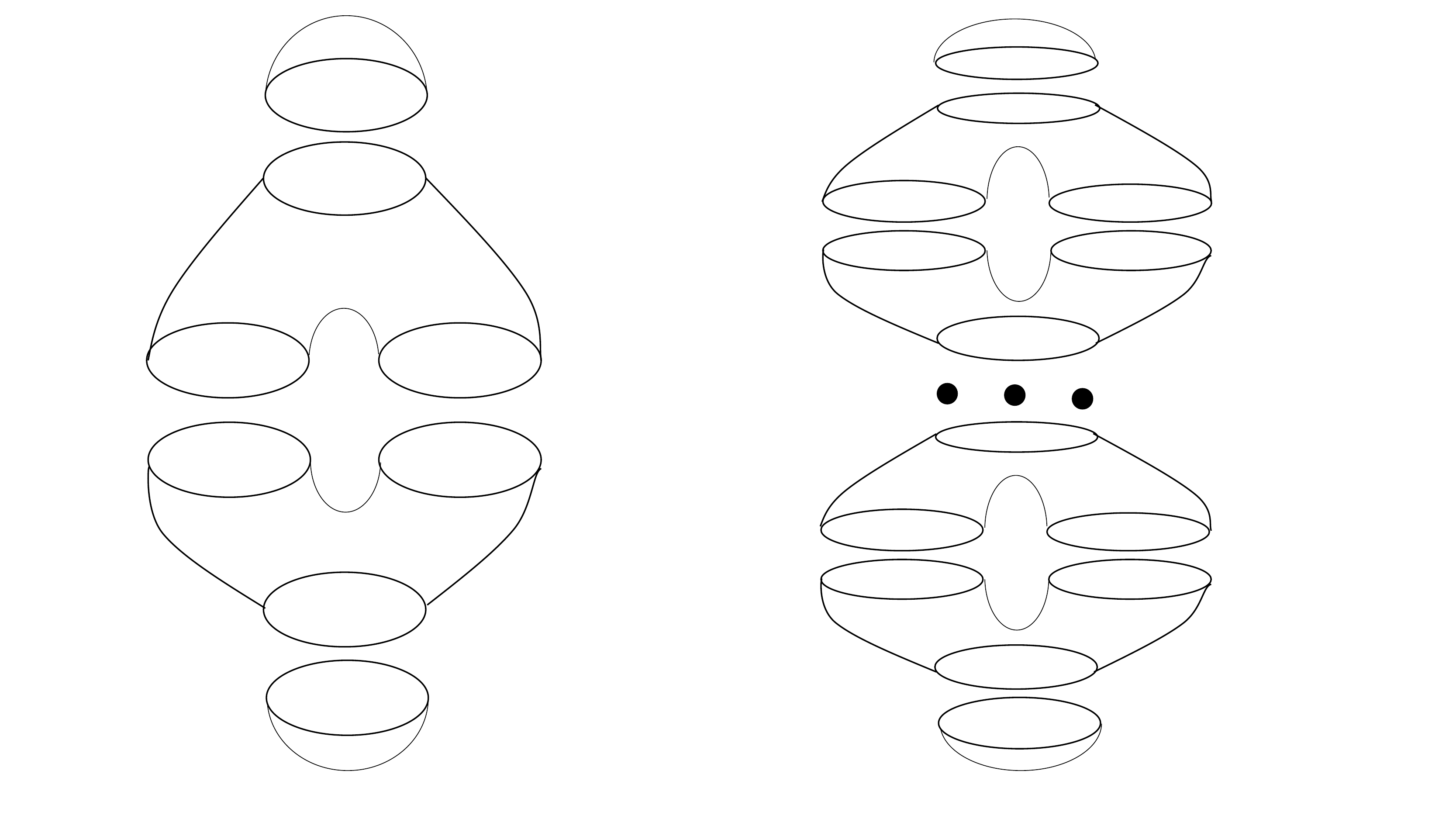}\\
  \caption{Cutting a Riemman surface into pairs of pants}\label{f1}
\end{figure}

In general, a Riemann surface with genus $g$ can be cut into $2g+2$ parts, which consist of a negative cap, $2g$ pairs of pant with a positive puncture and two symmetric negative puncture,  $2g$ pairs of pant with two positive puncture and a symmetric negative puncture and a positive cap. See figure \ref{f1}.

Let $\bar{D}$ be a Cauchy-Riemann operator on our Riemann surface $\Theta_{g, s^++s^-, t^++t^-}$ of genus $g$ with brake symmetry, which has $s^+$ positive puncture with brake symmetry, $s^-$ negative puncture with brake symmetry, $t^+$ pairs of positive punctures with only periodic condition, $t^-$ pairs of negative punctures with only periodic condition. We can glue $\Theta$ at each the punctures with opposite orientation caps of the same boundary condition. Then we get a closed Riemann surface. By the gluing formula, we have
\begin{align*}
&\text{index} \bar{D}|\Theta_{g, s^++s^-, t^++t^-} = \frac{n}{2}\chi(\Theta) + c_1(E)\\
&  - \sum_{i = 1}^{s^+}(\frac{n}{2} - \mu_1(q_i)) - \sum_{i' = 1}^{s^-}(\frac{n}{2} + \mu_1(q'_{i'})) - \sum_{j = 1}^{t^+}(n - \mu_{CZ}(p_j)) - \sum_{j' = 1}^{t^-} (n + \mu_{CZ}(p'_{j'}))\\
&= \frac{n}{2}(2-2g-s-s'-2t-2t')+ c_1(E)\\
& + \sum_{i = 1}^s\mu_1(q_i) -\sum_{i' = 1}^{s'}\mu_1(q'_{i'}) +\sum_{j = 1}^t \mu_{CZ}(p_j) - \sum_{j' = 1}^{t^-}\mu_{CZ}(p'_{j'})
\end{align*}

Each pair of periodic boundary are symmetric, we just count once. The Cauchy-Riemann operator $D$ on a positive cap $D^+$ of periodic boundary has index $n + \mu_{CZ}(S^\infty)$, this result comes from Schwarz\cite[\S 3.3.8 Proposition]{Schwarz}. And $D$ on negative cap $D^-$ has index $n - \mu_{CZ}(S^{-\infty})$.

\begin{figure}
  \centering
  \includegraphics[scale = 0.4]{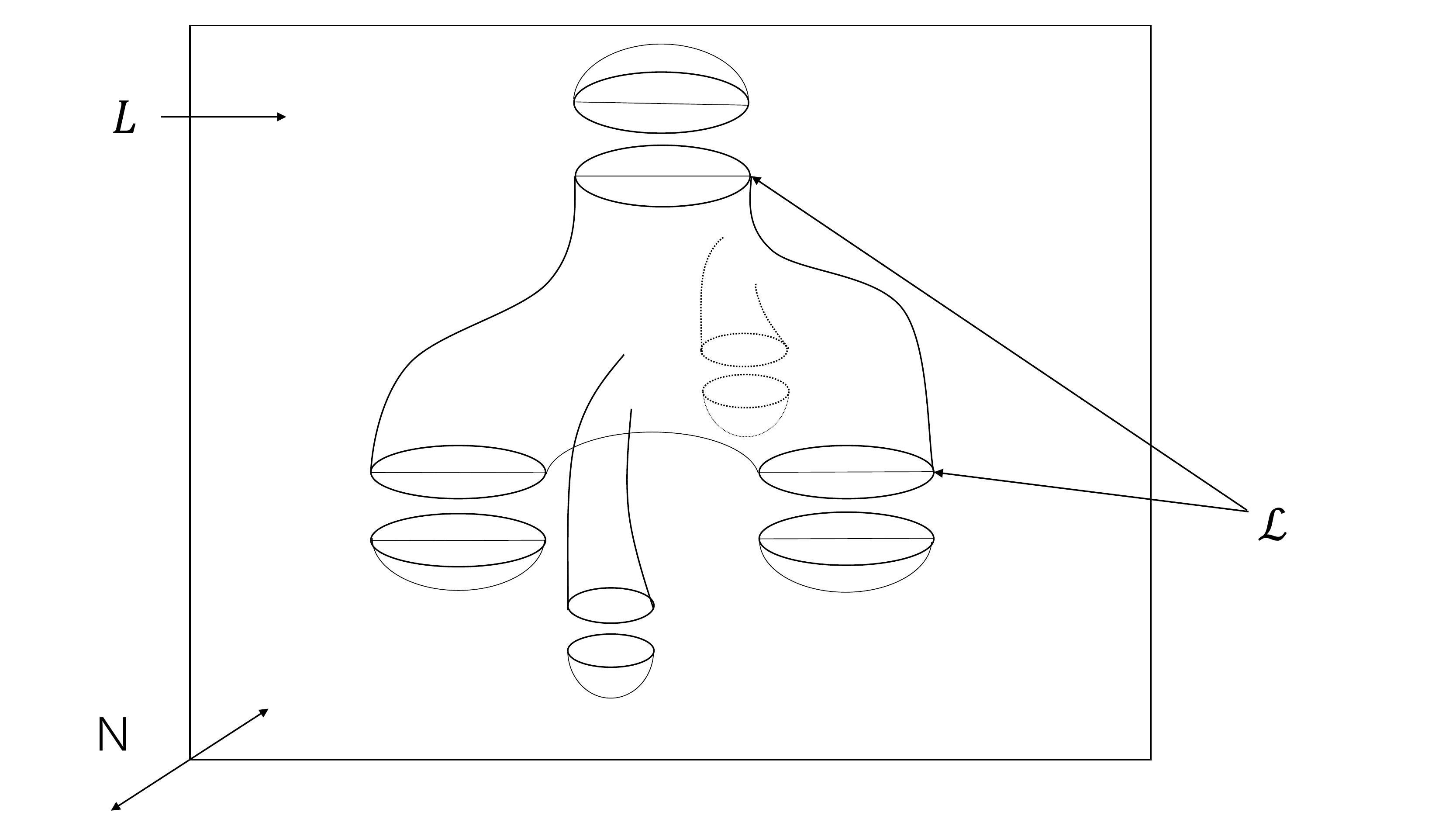}\\
  \caption{A pseudoholomorphic curve $F \in \mathcal{M}_{1+2,0+1}$}\label{f2}
\end{figure}

Here we give the picture of a pseudoholomorphic curve with brake symmetry, which has a positive symmetric puncture, $2$ negative symmetric punctures and a pair of negative punctures. See figure \ref{f2}.

Therefore we get our first main theorem, which is parallel to the result of \cite[\S 3.3.11 Theorem]{Schwarz}.
\begin{theo}\label{main1}
A Cauchy-Riemann operator $$\bar{D}: \bar{W}^{k, p}(\Theta, E) \to \bar{W}^{k-1, p}(\Theta, \Lambda^{0, 1} \otimes E)$$ has Fredholm index
\begin{align*}
& \frac{n}{2}(2- 2g - s^+-s^- -2t^+-2t^-)+ c_1(E) + \sum_{i = 1}^{s^+}\mu_1(q_i) -\sum_{i' = 1}^{s^-}\mu_1(q'_{i'})  \\
& +\sum_{j = 1}^{t^+}\mu_{CZ}(p_j) - \sum_{j' = 1}^{t^-}\mu_{CZ}(p'_{j'})
\end{align*}
\end{theo}

\begin{remark}
Taking a different trivialisation of $E$ over $\Theta$, the number $c_1$, $\mu_1(q_i), \mu_1(q'_{i'}), \mu_{CZ}(p_j), \mu_{CZ}(p'_{j'})$ can change but the number $c_1(E) + \sum_{i = 1}^s\mu_1(q_i) -\sum_{i' = 1}^{s'}\mu_1(q'_{i'}) +\sum_{j = 1}^t \mu_{CZ}(p_j) - \sum_{j' = 1}^{t'}\mu_{CZ}(p'_{j'})$ will not change.

We denote the first Chern number and Maslov index in a different trivialisation by $c_1(E'), \mu'_1, \mu'_{CZ}$. If we take a different trivialisation at a puncture, it is represented by a unitary loop $\phi(t)$ along the limit orbit. If $\phi(t)$ is a unitary loop and $\Phi$ is a symplectic path, then $\mu_{CZ}(\phi\Phi) = 2\mathrm{deg}(\phi) + \mu_{CZ}(\Phi)$, where $\mathrm{deg}(\phi) : \pi_1(\mathrm{U}(n)) \to \mathbf{Z}$ is induced by det$: \mathrm{U}(n) \to S^1$. See \cite[formula 3.48]{Schwarz}.

At a symmetric puncture because of our brake symmetry requirement for the trivialisation $T(\mathcal{N}\cdot)N = NT(\cdot)$, we need the condition $\phi(-t)N = N\phi(t)$. In appropriate coordinate such that $N$ is represented by the matrix $N_0$ on the fibre of complex bundle $E$. By Lemma \ref{4} below, we have $\mu_1(\phi\Phi) = \mathrm{deg}(\phi) + \mu_1(\Phi)$. It is well known that $c_1(E') - c_1(E) = -\sum_{*= i, i'}\mathrm{deg}(\phi_*)-\sum_{*=  j, j'}2\mathrm{deg}(\phi_*)$, which induces the equality
\begin{align*}
  & c_1(E) + \sum_{i = 1}^s\mu_1(q_i) -\sum_{i' = 1}^{s'}\mu_1(q'_{i'}) +\sum_{j = 1}^t \mu_{CZ}(p_j) - \sum_{j' = 1}^{t'}\mu_{CZ}(p'_{j'}) = \\
  & c_1(E') + \sum_{i = 1}^s\mu'_1(q_i) -\sum_{i' = 1}^{s'}\mu'_1(q'_{i'}) +\sum_{j = 1}^t \mu'_{CZ}(p_j) - \sum_{j' = 1}^{t'}\mu'_{CZ}(p'_{j'})
\end{align*}

\end{remark}

\begin{lemma}\label{4}
Let $({\bf R}^{2n}, \omega_0, J_0)$ be a linear space with symplectic form  $\omega_0 = \sum_{i= 1}^{n}dx_i\wedge dy_i$ and complex structure $J_0$. $({\bf R}^{2n}, J_0)$ can be seen as complex space ${\bf C}^{n}$ by isomorphism $z_i = (x_i, y_i)$, where $(x_1, \ldots, x_n, y_1, \ldots, y_n)$ and $(z_1, \ldots, z_n)$ are the coordinate of ${\bf R}^{2n}$ and ${\bf C}^{n}$. Given a unitary loop $\phi : S^1 \to \rm{U}(n)$ of period $\tau$ that satisfies $\phi(-t)N_0 = N_0\phi(t)$, we have $\mu_1(\phi\Phi) = \mathrm{deg}(\phi) + \mu_1(\Phi)$, where $\Phi : [0, \tau] \to \rm{Sp}(2n)$ is a symplectic path corresponding to a brake orbit of period $\tau$.
\end{lemma}

\begin{proof}
Firstly we have
\begin{equation*}
\phi(\frac{\tau}{2}+t) = N_0\phi(-\frac{\tau}{2}-t)N_0 = N_0\phi(\frac{\tau}{2}-t)N_0
\end{equation*}
the reason for the second equality is that the period of $\phi$ is $\tau$. So we have $\phi(\frac{\tau}{2}) = N_0 \phi(\frac{\tau}{2}) N_0$.

Let $\phi(\frac{\tau}{2}) = \left(
                                                                                                                  \begin{array}{cc}
                                                                                                                    A & B \\
                                                                                                                    C & D \\
                                                                                                                  \end{array}
                                                                                                                \right)
$. Then $N_0 \phi(\frac{\tau}{2}) N_0 = \left(
                                          \begin{array}{cc}
                                            A & -B \\
                                            -C & D \\
                                          \end{array}
                                        \right)
$. Hence we get $\phi(\frac{\tau}{2}) = \left(
                                                                                                                  \begin{array}{cc}
                                                                                                                    A & 0 \\
                                                                                                                    0 & D \\
                                                                                                                  \end{array}
                                                                                                                \right)$.

 $A, D$ are nondegenerate, i.e. $\text{det}(A), \text{det}(D) \neq 0$, since $\phi(t) \in \rm{U}(n)$.

Hence $\phi(\frac{\tau}{2})$ keeps the space $L_1= 0 \times {\bf R}^n$ and $L_2= {\bf R}^n \times 0$. By the same reason so does $\phi(0)$.

We can see $L_1 = 0 \times {\bf R}^n$ rotates in ${\bf C}^{n}$ by the action of $\phi(t)$. The rotation number of $\phi L_1$ equals $2\text{deg}(\phi)$.

From the definition of $\mu_1$, we have $\mu_1(\phi) =  \mu^{RS}(L_1, \phi[0, \frac{\tau}{2}] L_1)$.

Let $\phi_1(t) := N_0 \phi(\frac{\tau}{2}-t)N_0 : [0, \frac{\tau}{2}] \to \rm{U}(n)$, $\phi_2(t) := \phi(\frac{\tau}{2}-t)N_0 : [0, \frac{\tau}{2}] \to \rm{U}(n)$ and $\phi_3(t)= \phi(\frac{\tau}{2}-t) : [0, \frac{\tau}{2}] \to \rm{U}(n)$.

Using the equation  $\phi(\frac{\tau}{2} + t) = N_0 \phi(\frac{\tau}{2}-t)N_0$, $\mu^{RS}(L_1, \phi[\frac{\tau}{2}, \tau] L_1) = \mu^{RS}(L_1, \phi_1 L_1)$.

We have $\mu^{RS}(L_1, \phi_1 L_1) = -\mu^{RS}(N_0 L_1, \phi_2 L_1)$, since $N_0$ is antisymplectic $N^* \omega = -\omega$.

And $L_1$ is fixed by $N_0$, $N_0 L_1 = L_1$, we get $\mu^{RS}(N_0 L_1, \phi_2 L_1) = \mu^{RS}(L_1, \phi_3 L_1).$

According to the definition of $\mu^{RS}$, $\mu^{RS}(L_1, \phi_3 L_1) = -\mu^{RS}(L_1, \phi[0, \frac{\tau}{2}] L_1)$

By the results above, we conclude $\mu^{RS}(L_1, \phi[0, \tau] L_1) = \mu^{RS}(L_1, \phi[0, \frac{\tau}{2}] L_1) + \mu^{RS}(L_1, \phi[\frac{\tau}{2}, \tau] L_1) = 2\mu^{RS}(L_1, \phi[0, \frac{\tau}{2}] L_1)$.

From the definition of $\mu^{RS}$, $\mu^{RS}(L_1, \phi[0, \tau] L_1)$ equals rotation number of $\phi L_1$, which equals $2\text{deg}(\phi)$. This induces the equality $\mu_1(\phi) =  \mu^{RS}(L_1, \phi[0, \frac{\tau}{2}] L_1) = \text{deg}(\phi)$.

By the basic property of Maslov index, it is easy to get $\mu_1(\phi\Phi) = \mu_1(\phi) + \mu_1(\Phi)$, which yields $\mu_1(\phi\Phi) = \text{deg}(\phi) + \mu_1(\Phi)$.
\end{proof}

For our operator $\widetilde{D}_F $, $F \in \mathcal{M}_{s+s', t+t'}(q_1, \ldots, q_s; q'_1, \ldots, q'_{s'}; p_1, \ldots, p_t; p'_1,\\ \ldots, p'_{t'})$, we can use the result of theorem \ref{main1}. Here the vector bundle $E = F^*TM$ has two components: the $(a, \theta)$ component and the $\xi$ component. We have Morse-Bott construction, at each symmetric puncture the $(a, \theta)$ component has $\mu_1 = -\frac{1}{2}$ for positive punctures and $\mu_1 = \frac{1}{2}$ for negative punctures. Also at each symmetric pair of punctures the $(a, \theta)$ component has $\mu_{CZ} = -1$ for positive punctures and $\mu_{CZ} = 1$ for negative punctures. We know the $\xi$ component has index $\mu_1(q)$ at each symmetric puncture $q$ and index $\mu_{CZ}(p)$ at each symmetric pair of punctures $p, Np(-t)$. Therefore we have $\mu_1(q_i)|F^*TM = \mu_1(q_i) - \frac{1}{2}$, for $i = 1, \ldots, s$, $\mu_1(q'_{i'})|F^*TM = \mu_1(q'_{i'}) + \frac{1}{2}$,for $i = 1, \ldots, s'$, $\mu_{CZ}(p_j)|F^*TM = \mu_{CZ}(p_j) - 1$, for $j = 1, \ldots, t$, $\mu_{CZ}(p'_{j'})|F^*TM = \mu_{CZ}(p'_{j'}) + 1$, for $j = 1, \ldots, t$,$\mu_1(q_i), \mu_1(q'_{i'}), \mu_{CZ}(p_j), \mu_{CZ}(p'_{i'})$ are our original definition of indices for the orbits $q_1, \ldots, q_s; q'_1, \ldots, q'_{s'}; p_1, \ldots, p_t; p'_1, \ldots, p'_{t'}$.

Therefore we get index $\widetilde{D}_F = \frac{n}{2}(2-2g-s-s'-2t-2t')+ \sum_{i = 1}^s\mu_1(q_i) -\sum_{i' = 1}^{s'}\mu_1(q'_{i'}) +\sum_{j = 1}^t \mu_{CZ}(p_j) - \sum_{j' = 1}^{t'} \mu_{CZ}(p'_{j'})-\frac{1}{2}N$, $N = s+s'+ 2t+2t'$.

And for the total operator $D'_F = \pi \circ d\bar{\partial}_F$, $D'_F$ has domain ${\bf R}^{s+s'+2t+2t'} \oplus \bar{W}^{d,k,p}(\Theta, F^*TM)$. We have calculated the index $D_F = D'_F|\bar{W}^{d,k,p}$. A elementary property of Fredholm operator is index $S \circ T = \text{index}\ S + \text{index}\ T$, where $T : U \to V, S : V \to W$ are Fredholm operators. $D'_F|\bar{W}^{d,k,p} = \widetilde{D}_F \circ i$, where $i$ is the embedding of $\bar{W}^{d,k,p}$ into the second component of the range $i : \bar{W}^{d,k,p} \to {\bf R}^{s+s'+2t+2t'} \oplus \bar{W}^{d,k,p}(\Theta, F^*TM), i(u) = (0, u)$. It is easy to check $i$ is a Fredholm operator. We can get index $D'_F = \text{index}\ \widetilde{D}_F + N =\frac{n}{2}(2-2g-s-s'-2t-2t')+ \sum_{i = 1}^s\mu_1(q_i) -\sum_{i' = 1}^{s'}\mu_1(q'_{i'}) +\sum_{j = 1}^t \mu_{CZ}(p_j) - \sum_{j' = 1}^{t'} \mu_{CZ}(p'_{j'}) + \frac{1}{2}N$.

We want to introduce some facts of Riemann surface. Let $\overline{\Theta}$ be a Riemann surface. Let $\mathcal{S}(\overline{\Theta})$ be the space of all complex structure of $\overline{\Theta}$, $\textrm{Diff}^+_0(\overline{\Theta})$ be the group of orientation-preserving diffeomorphisms of $\overline{\Theta}$ that are isotopic to the identity. The Teichm\"{u}ller space $\mathcal{T}(\overline{\Theta})$ for Riemann surface $\overline{\Theta}$ is the quotient space $\mathcal{S}(\overline{\Theta})/\textrm{Diff}^+_0(\overline{\Theta})$. See \cite{Teichmuller}, \cite{Abikoff}, \cite{Hubbard}.

Recall our definition of a pair $(\Theta, \mathcal{N})$, $\Theta = \Xi_g \backslash \{u_1, \ldots, u_{s^+}; u'_1, \ldots, u'_{s^-}; v_1, \ldots,\\ v_{t^+}, w_1, \ldots, w_{t^+}; v'_1, \ldots, v'_{t^-}, w'_1, \ldots, w'_{t^-}\}$ with $s+s'$ punctures $ u_1, \ldots, u_{s^+}, u'_1, \ldots,\\ u'_{s^-}$ such that $\mathcal{N}u_i = u_i, \mathcal{N}u'_{i'} = u'_{i'}, i= 1, \ldots, s^+, i'= 1, \ldots, s^-$, and $t^+ +t^-$ pairs of punctures $v_1, \ldots, v_{t^+}, w_1, \ldots, w_{t^+}, v'_1, \ldots, v'_{t^-}, w'_1, \ldots, w'_{t^-}$ such that $\mathcal{N}v_k = w_k, \mathcal{N}v_{k'} = w_{k'}$, for $k = 1, \ldots, t^+, k' = 1, \ldots, t^-$ and with the symmetry $\mathcal{N}j = -j\mathcal{N}$. $\mathcal{L}$ is the fixed set of $\mathcal{N}$, which is nonempty. And $\Xi_g \backslash \mathcal{L}$ has two connected components. We denote the two connected component of $\Xi_g \backslash \mathcal{L}$ by $\Xi'_g$ and $\Xi''_g$. Without loss of generality, we assume that the punctures $v_1, \ldots, v_{t^+}, v'_1, \ldots, v'_{t^-}$ are in $\Xi'_g$ and $w_1, \ldots, w_{t^+}, w'_1, \ldots, w'_{t^-}$ are in $\Xi''_g$. We define $\Theta' := \Xi'_g \bigcup \mathcal{L} \backslash \{u_1, \ldots, u_{s^+}; u'_1, \ldots, u'_{s^-}; v_1, \ldots, v_{t^+}; v'_1, \ldots, v'_{t^-} \}$ and $\Theta'' := \Xi''_g \bigcup \mathcal{L} \backslash \{u_1, \ldots, u_{s^+}; u'_1, \ldots, u'_{s^-}; w_1, \ldots, w_{t^+}; w'_1, \ldots, w'_{t^-} \}$, which are Riemann surfaces with boundary and have punctures on the boundary and in the interior.

The complex structure space with brake symmetry $\mathcal{N}j = -j\mathcal{N}$ on the space $(\Theta, \mathcal{N})$ is equivalent to the complex structure space on the space $\Theta'$. On one side, by restriction $j$ from $\Theta$ to $\Theta'$ we can get a complex structure on $\Theta'$. On the other side, from a complex structure on $\Theta'$, we can extend the complex structure by conjugation to $\Theta''$. Then gluing $\Theta'$ and $\Theta''$ along the boundary, we will get a complex structure on the space $\Theta$. For the details, see Abikoff \cite[page 44]{Abikoff}. We denote the Teichm\"{u}ller space with brake symmetry on $\Theta_{g, s^+ +s^-,t^+ +t^-}$ by $\mathcal{T}_{g, s^+ +s^-,t^+ +t^-}$.

An important fact is,
\begin{lemma}
	If $3g + s^+ +s^- + 2t^+ + 2t^- -3 > 0$, the real dimension of the Teichm\"{u}ller space $\mathcal{T}_{g, s^+ +s^-,t^+ +t^-}$ is $3g + s^+ + s^- +2t^+ +2t^- -3$. If $3g + s^+ + s^- +2t^+ +2t^- -3 \leq 0$, the real dimension of the Teichm\"{u}ller space $\mathcal{T}_{g, s^+ +s^-,t^+ +t^-}$ is $0$.
\end{lemma}

\begin{proof}
	According to the standard Teichm\"{u}ller theory\cite{Abikoff}, \cite{Hubbard}, the real dimension of the Teichm\"{u}ller space ${\mathcal T}_{g, k}$ for a Riemann surface $\Theta_{g, k}$ of genus $g$ with $k$ punctures is $6g + 2(k-3)$, when $6g + 2(k-3) > 0$; and it is $0$, when $6g + 2(k-3) \leq 0$. The Teichm\"{u}ller space ${\mathcal T}_{g, k}$ is isomorphic to the space $Q(\Omega^2)$, where $Q(\Omega^2)$ denotes the space of meromorphic quadratic differentials, such that at each of the punctures, the meromorphic quadratic differentials has at worst a first order pole.
	
	In our case the complex structure in $\Theta$ has the symmetry $\mathcal{N}j = -j\mathcal{N}$. The correspondent meromorphic quadratic differentials $\omega \in Q(\Omega_\Theta^2)$ has the symmetry $\mathcal{N}^*\omega = \omega$, then it is real on the fixed set $\mathcal{L}$. We define $Q_1= \{\omega|\mathcal{N}^*\omega = \omega\}$ and $Q_2= \{\omega|\mathcal{N}^*\omega = -\omega\}$, $Q_1$ is real on $\mathcal{L}$ and $Q_2$ is imaginary on $\mathcal{L}$. $Q_1$ and $Q_2$ are isomorphic by multiplication of complex number $j$. Every meromorphic quadratic differential $\omega$ can be written as $\omega = \omega_1 + \omega_2, \omega_1\in Q_1, \omega_2 \in Q_2$ by let $\omega_1 = \frac{\omega(\cdot) + N\omega(\mathcal{N} \cdot)}{2}, \omega_2 = \frac{\omega(\cdot) - N\omega(\mathcal{N} \cdot)}{2}$. So dim $Q_1 = \frac{1}{2} \text{dim} Q =3g + k-3$, where $k =  s^+ + s^- +2t^+ +2t^-$.
	
Hence dim $\mathcal{T}_{g, s^+ +s^-,t^+ +t^-} = \frac{1}{2}$ dim $\mathcal{T}_{g, s^+ + s^- +2t^+ +2t^-} =\\
	\begin{cases}
	3g + s^+ +s^-,t^+ +t^-- 3,  & if\ 3g + s^+ + s^- +2t^+ +2t^- - 3 > 0,\\
	0, & if\ 3g + s^+ + s^- +2t^+ +2t^- - 3 \leq 0.
	\end{cases}$
\end{proof}

We define $\textrm{Aut}(\Theta_{g, s^+ +s^-,t^+ +t^-})$ to be the automorphism group of\\ $(\Theta_{g, s^+ +s^-,t^+ +t^-}, j)$, where $j$ is a complex structure having brake symmetry.
Another fact is
\begin{lemma}
	If $3g + s^+ + s^- +2t^+ +2t^- -3 > 0$, the real dimension of Automorphism group $\textrm{Aut}(\Theta_{g, s^+ +s^-,t^+ +t^-})$ is $0$. If $3g + s^+ + s^- +2t^+ +2t^- -3 \leq 0$, the real dimension of $\textrm{Aut}(\Theta_{g, s^+ +s^-,t^+ +t^-})$ is $3- 3g - s^+ - s^- -2t^+ -2t^-$.
\end{lemma}

\begin{proof}
If $g \geq 1$, then $3g + s^+ + s^- +2t^+ +2t^- -3 \geq 0$, we know the Automorphism group should be either trivial or has dimension $0$.

Therefore we only have to consider the case $g=0$.
We have $\textrm{Aut}(\Theta_{0, s^+ +s^-,t^+ +t^-})\\ = \textrm{Aut}(\Theta'_{0, s^+ +s^-,t^+ +t^-})$, where $\Theta'_{0, s^+ +s^-,t^+ +t^-}$ is holomorphic to the unit disk in $\bf C$ with $s^+ +s^-$ punctures on the boundary and $t^+ +t^-$ punctures in the interior.
	
	The automorphism group of the unit disc in $\bf C$ has real dimension $3$. Each time we add one puncture on the boundary, the automorphism group will lose one dimension to keep the puncture fixed; and each time we add one puncture inside the disk, the automorphism group will lose two dimensions to keep the puncture fixed. Eventually, when $s^+ + s^- +2t^+ +2t^- - 3 > 0$, the automorphism group will be trivial.
\end{proof}

If $3g + s+s'+2t+2t'-3 > 0$, the dimension of  $\mathcal{T}_{g, s+s', t+t'} = 3g + s+s'+2t+2t' - 3$, the dimension of $Aut(\Theta_{g, s+s', t+t'})$ is $0$; if $3g + s+s'+2t+2t'-3 \leq 0$, the dim $\mathcal{T}_{g, s+s', t+t'} = 0$, the dimension of $Aut(\Theta_{g, s+s', t+t'})$ is $3-3g-( s+s'+2t+2t')$.

Hence, index of $D'_F$ $+$ dimension of Teichm\"{u}ller space $\mathcal{T}_{g, s+s', t+t'}$ $-$ dimension of automorphism group of Riemann surface with punctures $\Theta_{g, s+s', t+t'} =$
\begin{equation*}
\frac{n-3}{2}(2-2g-s-s'-2t-2t')+ \sum_{i = 1}^s\mu_1(q_i) -\sum_{i' = 1}^{s'}\mu_1(q'_{i'}) +\sum_{j = 1}^t \mu_{CZ}(p_j) - \sum_{j' = 1}^{t'} \mu_{CZ}(p'_{j'})
\end{equation*}

We define the virtual dimension of $\mathcal{M}_{s+s', t+t'}(q_1, \ldots, q_s; q'_1, \ldots, q'_{s'}; p_1, \ldots, \\p_t; p'_1, \ldots, p'_{t'})$ as index $D'_F +$ dim $\mathcal{T}_{s+s', t+t'}-$ dim Auto$(\Theta) =$
\begin{equation*}
\frac{n-3}{2}(2-2g-s-s'-2t-2t')+ \sum_{i = 1}^s\mu_1(q_i) -\sum_{i' = 1}^{s'}\mu_1(q'_{i'}) +\sum_{j = 1}^t \mu_{CZ}(p_j) - \sum_{j' = 1}^{t'} \mu_{CZ}(p'_{j'})
\end{equation*}

\begin{theo}[Main Theorem]\label{e}
	The virtual dimension of moduli space

$\mathcal{M}_{g, s+s', t+t'}(q_1, \ldots, q_s; q'_1, \ldots, q'_{s'}; p_1, \ldots, p_t; p'_1, \ldots, p'_{t'})$ is
	$$\frac{n-3}{2}(2-2g-s-s'-2t-2t')+ \sum_{i = 1}^s\mu_1(q_i) -\sum_{i' = 1}^{s'}\mu_1(q'_{i'}) +\sum_{j = 1}^t \mu_{CZ}(p_j) - \sum_{j' = 1}^{t'} \mu_{CZ}(p'_{j'})$$.
\end{theo}

Especially, when we take $g=0, s=1, t=0$, the virtual dimension of $\mathcal{M}_{1+s', t'}(q; q'_1, \ldots, q'_{s'};\emptyset; p'_1, \ldots, p'_{t'})$ is
\begin{equation*}
\frac{n-3}{2} + \mu_1(q) - \sum_{i=1}^{s}\left( \frac{n-3}{2} + \mu_1(q'_{i'})\right) -\sum_{j'=1}^{t'}\left((n-3) + \mu_{CZ}(q'_{j'})\right)
\end{equation*}

So we define the degree of a brake orbit $q$ by $\mu_1(q) + \frac{n-3}{2}$, and the degree of a closed Reeb orbit $p$ by $\mu_{CZ}(p) + n-3$.

We have the virtual dimension of  $\mathcal{M}_{1+s', t'}(q; q'_1, \ldots, q'_{s'};\emptyset; p'_1, \ldots, p'_{t'}) = |q| - \sum_{{i'}=1}^{s'}|q'_{i'}| - \sum_{j'=1}^{t'}|p'_{j'}|$

\begin{theo}
	The virtual dimension of moduli space $\mathcal{M}_{1+s', t'}(q; q'_1, \ldots, q'_{s'};\\ \emptyset; p'_1, \ldots, p'_{t'})$ is
	$$\frac{n-3}{2} + \mu_1(q) - \sum_{i=1}^{s}\left( \frac{n-3}{2} + \mu_1(q'_{i'})\right) -\sum_{j'=1}^{t'}\left((n-3) + \mu_{CZ}(q'_{j'})\right)$$
\end{theo}
In this theorem $q;q'_1, \ldots, q'_{s'}$ be brake orbit and  $p'_1, \ldots, p'_{t'}$ be closed Reeb orbit. In order to define contact homology for brake orbits, we must take care of the phenomena of  bubbling off. The limit of a series of the pseudoholomorphic curve with brake symmetry, which have only brake orbits boundary, can bubble off to a pseudoholomorphic curve with brake symmetry which has negative pairs of periodic orbits boundary.

\section{Outlook}
Similar to the construction of contact homology, we want to use the same procedure to define the contact homology for brake orbits by counting pseudoholomorphic curves with brake symmetry in ${\mathcal M}_{1+s', 0+0}$. We call this new homology real contact homology. At the same time, we can also define Real cylindrical contact homology by counting pseudoholomorphic cylinders with brake symmetry.

Contact homology has two versions: equivariant and nonequivariant contact homology. Nonequivariant contact homology uses $S^1$ dependent almost complex structures $J(t)$. When there is a filling for the contact manifold $(\Sigma, \alpha)$, i.e. there exists a symplectic manifold $(M, \omega)$, such that $\partial M = \Sigma$. At the same time, there exists a vector field $X$ in a neighborhood of $\Sigma$, such that $L_{X}\omega = \omega, \alpha = \omega(X, \cdot)$. The equivariant(nonequivariant) contact homology of $\Sigma$ is isomorphic to the equivariant(nonequivariant)symplectic homology of $M$\cite{BO}. In a future project for which the index computations in this paper are crucial is to define nonequivariant real contact homology as well.

J.Kim, S.Kim and Kwon have already defined equivariant wrapped Floer homology\cite{Kim}. We guess it should be equivalent to the equivariant real contact homology and we hope to prove this statement in the future.

In addition, because we know that a brake orbit is equivalent to a path satisfies the equations (\ref{104})-(\ref{105}), we can define a homology by counting pseudoholomorphic strips with such paths as boundary. And there are also two versions, equivariant and nonequivariant. We call this homology equivariant(nonequivariant) strip-like real contact homology.

More ambitiously, we can consider analogy of the embedded contact homology and symplectic field theory with brake symmetry. We call the corresponding theory real embedded contact homology and real symplectic field theory, which we abbreviate by real ECH and real SFT.

All these theories are interesting topics for future research.

\section*{Acknowledgements}
This paper is accomplished when the first author is visiting University of Augsburg, Germany. The first author is supported by China Scholarship Counsil(CSC) No.201806200130, LPMC of Ministry of Education of China, Nankai Zhide Foundation and Nankai University. The second author is supported by NSFC Grants 11971245 and 11771331, LPMC of Ministry of Education of China, Nankai Zhide Foundation and Nankai University. Professor Frauenfelder has given many valuable advises to this paper. The first author also thanks professor Frauenfelder's hospitality and fruitful discussions when he stays in Germany.


\bibliography{paper10}

\end{document}